\documentclass[11pt]{article}
\usepackage[tbtags]{amsmath}
\usepackage{amssymb}
\usepackage{amsthm}
\usepackage[misc]{ifsym}
\usepackage{cases}
\usepackage{mathrsfs}
\usepackage{color}
\usepackage{graphicx}
\usepackage{subfigure}
\usepackage{xcolor}
\usepackage{tikz}
\usetikzlibrary{arrows,shapes,chains}
\numberwithin{equation}{section}
\normalsize
\title{\bf Closed-Loop Solvability of Stochastic Linear-Quadratic Optimal Control Problems with Poisson Jumps
\thanks{This work is supported by National Natural Science Foundations of China (Grant Nos. 11971266, 11831010, 11571205), and Shandong Provincial Natural Science Foundations (Grant Nos. ZR2020ZD24, ZR2019ZD42).}}
\author{\normalsize Zixuan Li\thanks{\it School of Mathematics, Shandong University, Jinan 250100, P.R. China, E-mail: 201812064@mail.sdu.edu.cn},\quad Jingtao Shi\thanks{\it Corresponding author. School of Mathematics, Shandong University, Jinan 250100, P.R. China, E-mail: shijingtao@sdu.edu.cn}}
\newtheorem{mypro}{Proposition}[section]
\newtheorem{mythm}{Theorem}[section]
\newtheorem{mydef}{Definition}[section]

\begin{document}
\maketitle

\noindent{\bf Abstract:}\quad This paper is concerned with the stochastic linear-quadratic optimal control problem with Poisson jumps. The coefficients in the state equation and the weighting matrices in the cost functional are all deterministic but are allowed indefinite. The notion of closed-loop strategies is introduced, and the optimal closed-loop strategy is characterized by a Riccati integral-differential equation and a backward stochastic differential equation with Poisson jumps.

\vspace{2mm}

\noindent{\bf Keywords:}\quad stochastic linear-quadratic optimal control, Poisson random measure, backward stochastic differential equation with Poisson jumps, Riccati integral-differential equation, closed-loop solvability

\vspace{2mm}

\noindent{\bf Mathematics Subject Classification:}\quad 49N10, 49K45, 60H10, 93E20

\section{Introduction}

\par The {\it stochastic linear-quadratic} (SLQ) optimal control problem plays an extremely important role in modern control theory and methodology, for its elegant structure of solutions and wide applications in engineering, finance and network, etc. More importantly, many nonlinear stochastic control problems can be reasonably approximated by the SLQ problems. Literatures for SLQ optimal control problems, refer to Wonham \cite{Wonham1968}, Bismut \cite{Bismut1976}, Bensoussan \cite{Bensoussan1981}, Peng \cite{Peng1992}, Chen et al. \cite{ChenLiZhou1998}, Chen and Zhou \cite{ChenZhou2000}, Chen and Yong \cite{ChenYong2001}, Ait Rami et al. \cite{AMZ2001}, Tang \cite{Tang2003}, Yu \cite{Yu2013}, Tang \cite{Tang2015}, Sun et al. \cite{SLY2016}, Sun and Yong \cite{SY2018}, Sun et al. \cite{SunXiongYong2021}, etc. for journal papers and Davis \cite{Davis1977}, Anderson and Moore \cite{AndersonMoore1989}, Yong and Zhou \cite{YongZhou1999}, Sun and Yong \cite{SunYong2020} for monographs.

\par In the above works, stochastic systems are modelled by Browonian motions. However, in reality, Brownian noises are usually inadequate in mathematical modeling sense. For example, it is particularly appropriate to use stochastic systems with Poisson jumps or L\'{e}vy jumps to describe the large fluctuations in the stock market (Merton \cite{Merton1976}, Kou \cite{Kou2002}, Cont and Tankov \cite{ContTankov2004}, \O ksendal and Sulem \cite{OS2005}, Lim \cite{Lim2005}, Hanson \cite{Hanson2007}). Moreover, from a mathematical point of view, there exist essential differences between stochastic systems with and without jumps.

{\it SLQ optimal control problems with Poisson jumps} (SLQP optimal control problems) are also researched by many authors. Tang and Hou \cite{TangHou2002} studied an SLQP optimal control problem with a free Poisson process and a Gaussian white-noised observation, the optimal control is explicitly solved using the partially observed maximum principle. Wu and Wang \cite{WuWang2003} studied a kind of SLQP optimal control problems, the explicit form of optimal controls is obtained by the solutions to a {\it forward-backward stochastic differential equation with Poisson jumps} (FBSDEP) and a generalized Riccati equation system. Hu and \O ksendal \cite{HuOksendal2008} studied an SLQP optimal control problem with partial information. Meng \cite{Meng2014} considered an SLQP optimal control problem with random coefficients. The state feedback representation was obtained for the open-loop optimal control by a matrix-valued {\it backward stochastic Riccati equations with jumps} (BSREJ), and the solvability of it in a special case was discussed. The solvability of BSREJ in general case was studied in Zhang et al. \cite{ZhangDongMeng2020}. Note that Li et al. \cite{LiWuYu2018} used a so-called relax compensator to describe indefinite BSREJ and investigated the solvability of BSREJ in some special cases. Moon and Chung \cite{MoonChung2021} studied the indefinite SLQP optimal control problem with random coefficients by a completion of squares approach.

\par Our interest in this paper is the closed-loop solvability of the SLQP optimal control problem, which as far as we know, is not researched in the literature. In 2014, the notions of open-loop and closed-loop solvabilities for SLQ optimal control problems are introduced in Sun and Yong \cite{SY2014}, which is a special case when only one player/controller is considered in the LQ zero-sum stochastic differential game. Sun et al. \cite{SLY2016} further gave more detailed characterizations of the closed-loop solvability for SLQ optimal control problems. Sun and Yong \cite{SY2018} obtained the equivalence of open-loop and closed-loop solvabilities for the SLQ optimal control problem in an infinite horizon. See also their paragraph \cite{SunYong2020}. Li et al. \cite{LSY2016} studied the existence of an optimal closed-loop strategy for an SLQ optimal control problem of mean-field's type. Lv \cite{Lv2019,Lv2020} researched the closed-loop solvabilities of SLQ optimal control problems for systems governed by {\it stochastic evolution equations} (SEEs) and SEEs of mean-field's type, respectively. Tang et al. \cite{TangLiHuang2020} studied the open-loop and closed-loop solvability for indefinite SLQP optimal control problem of mean-field's type. Sun et al. \cite{SunXiongYong2021} considered the indefinite SLQ optimal control problem with random coefficients, and investigated the closed-loop representation of open-loop optimal controls.

\par Our work distinguishes with the existing results in the following respects. (1) We consider an SLQP optimal control problem with deterministic coefficients in a general framework ({\bf Problem (SLQP)} in Section 2), where the weighting matrices in the cost functional are allowed to be indefinite. Moreover, cross-product terms in the control and state processes are present in the cost functional. Non-homogenous terms are also appear in the controlled state equation and cost functional. The model considered in this paper is a nontrivial generalization of those in \cite{SY2014,SLY2016}. (2) Characterization of the closed-loop solvability for the SLQP optimal control problem is obtained, via some {\it Riccati integral-differential equation} (RIDE). For the SLQ optimal control problem without Poisson jumps, Sun et al. \cite{SLY2016} first found two matrix-valued SDE of $\mathbb{X}(\cdot)$, $\mathbb{Y}(\cdot)$, then they applied It\^o's formula to find $\mathbb{X}^{-1}(\cdot)$. The solution to the related Riccati equation is defined as $P(\cdot)=\mathbb{Y}(\cdot)\mathbb{X}^{-1}(\cdot)$. But this method fails in our SLQP optimal control Problem, since the Poisson jumps appear in the controlled system and the difficulty is encountered. In detail, when we take the inverse of the matrix-valued SDEP as that of the matrix-valued SDE in \cite{SLY2016} and applying It\^{o}--Wentzell's formula (\O ksendal and Zhang \cite{OksendalZhang2007}), term like $\big[F(e)+G(e)\bar{\Theta}+1\big]^{-1}$ will appear in $\mathbb{X}^{-1}(\cdot)$. But since we don't have any restrictions on the coefficients in our system and $\bar{\Theta}(\cdot)$ is the closed-loop optimal strategy which we are going to seek, there is no reason to presume arbitrarily that $F(e)+G(e)\bar{\Theta}+1$ is invertible. Thank for Lv \cite{Lv2019,Lv2020}, we overcome this difficulty by transforming the original Problem (SLQP) into a problem of solving the open-loop optimal control of {\bf Problem (SLQP)$^0_{\bar{\Theta}}$} in Section 2. Thus, a Lyapunov integral-differential equation (\ref{Lyapunov}) is given first and then the RIDE (\ref{Riccati}) is obtained. Note that the technique used in this paper is also different from that in Tang et al. \cite{TangLiHuang2020}, where a matrix minimum principle by Athans \cite{Athans1968} was used when dealing with the SLQP problem of mean-field's type.

\par The rest of this paper is organized as follows. Section 2 gives some preliminaries and problem formulation. Section 3, main results are presented and proved in detail. Finally, in Section 4 some concluding remarks are given.

\section{Problem formulation and preliminaries}

%Problem(SLQP)
First, let us introduce some notations which will be used throughout this paper.

Let $T>0$ be a finite time duration. Let $\mathbb{R}^{n \times m}$ be the set of all $(n \times m)$ matrices, and $\mathbb{S}^n$ be the set of all $(n \times n)$ symmetric matrices. We let $I$ be the identity matrix with a suitable size. We use $\langle\cdot\,,\cdot\rangle$ for inner products in possibly different Hilbert spaces, and denote by $|\cdot|$ the norm induced by $\langle\cdot\,,\cdot\rangle$. Let $M^\top$ and $\mathcal{R}(M)$ be the transpose and range of a matrix $M$, respectively. For $M,N\in\mathbb{S}^n$, we write $M\geqslant N$ (respectively, $M>N$) for $M-N$ being positive semi-definite (respectively, positive definite). Let $M^\dag$ denote the pseudo-inverse of a matrix $M\in\mathbb{R}^{m\times n}$, which is equal to the inverse $M^{-1}$ of $M\in\mathbb{R}^{n\times n}$ if it exists. See Penrose \cite{Penrose1955} for some basic properties of the pseudo-inverse.

For any Banach space $H$ (for example, $H=\mathbb{R}^n,\mathbb{R}^{n \times m},\mathbb{S}^n$) and $t \in [0,T)$, let $L^p(t,T;H)\,(1 \leqslant p \leqslant \infty)$ be the space of all $H$-valued functions that are $L^p$-integrable on $[t,T]$, $C([t,T];H)$ be the space of all $H$-valued continuous functions on $[t,T]$, and $L^\infty(t,T;H)$ be the space of Lebesgue measurable, essentially bounded functions from $[t,T]$ into $H$.

Let $(\Omega,\mathcal{F},\mathbb{F},\mathbb{P})$ be a completed filtered probability space, where $\mathbb{F}=\{\mathcal{F}_t\}_{t\geqslant 0}$ is filtration generated by the following two mutually independent stochastic processes and augmented by all the $\mathbb{P}$-null sets in $\mathcal{F}$:
\begin{itemize}
	\item A standard one-dimensional Brownian motion $W=\{W(t);0 \leqslant t < \infty \}$.
	\item A Poisson random measure $\mathnormal{N}$ on $\mathbf{E} \times \mathbf{R}_+$, where $\mathbf{E}  \subset \mathbf{R}^l $ is a nonempty open set equipped with its Borel field $\mathcal{B}(\mathbf{E})$, with the compensator $\hat{\mathnormal{N}}(dedt)=\pi(de)dt$, satisfying $\pi(A) \textless \infty$, such that $\tilde{\mathnormal{N}}(A\times[0,t])=(\mathnormal{N}-\hat{\mathnormal{N}})(A\times[0,t])_{t \geqslant 0}$ is a martingale for all $A \in \mathcal{B}(\mathbf{E})$. $\pi$ is assumed to be a $\sigma$-finite measure on $(\mathbf{E},\mathcal{B}(\mathbf{E}))$ and is called the characteristic measure.
\end{itemize}
For $t \in [0,T)$, we denote
\begin{equation*}
\begin{aligned}
&L^2_{\mathcal{F}_s}(\Omega;H) =\Big\{\xi:\Omega \to H\,\big|\,\xi\,\,\mbox{is}\,\,\mathcal{F}_s\mbox{-measurable},\,\,\,\mathbb{E}|\xi|^2<\infty \Big\},\,\,\, s\in[t,T],\\
&L^2_{\mathbb{F}}(t,T;H) =\bigg\{f:[t,T]\times\Omega\to H\,\big|\,f\,\,\mbox{is}\,\,\mathbb{F}\mbox{-progressively measurable},\,\,\mathbb{E}\int_t^T|f(s)|^2ds<\infty \bigg\},\\
&L^2_{\mathbb{F},p}(t,T;H) =\bigg\{f:[t,T]\times\Omega\to H\,\big|\,f\,\,\mbox{is}\,\,\mathbb{F}\mbox{-predictable},\,\,\mathbb{E}\int_t^T|f(s)|^2ds<\infty \bigg\},\\
&F^\infty_{p}(t,T;H) =\bigg\{f:[t,T]\times\mathbf{E}\to H\,\big|\,f\,\,\mbox{satisfies}\,\, \mathop{\sup}\limits_{t\leqslant s\leqslant T,\,e\in\mathbf{E}}|f(s,e)|<\infty \bigg\},\\
&F^2_{\mathbb{F},p}(t,T;H) =\bigg\{f:\Omega\times[t,T]\times\mathbf{E}\to H\,\big|\,f\,\,\mbox{is}\,\,\mathbb{F}\mbox{-predictable},
\mathbb{E}\int_t^T\int_\mathbf{E}|f(\cdot,s,e)|^2\pi(de)ds<\infty \bigg\}.
\end{aligned}
\end{equation*}

We consider the following controlled linear SDEP on $[t,T]$:
\begin{equation}\label{state}\left\{
\begin{aligned}
dX(s)&=\big[A(s)X(s)+B(s)u(s)+b(s)\big]ds+\big[C(s)X(s)+D(s)u(s)+\sigma(s)\big]dW(s)\\
     &\qquad+\int_\mathbf{E}\big[F(s,e)X(s-)+G(s,e)u(s)+f(s,e)\big]\tilde{\mathnormal{N}}(deds),\quad s \in [t,T],\\
X(t)&=x,
\end{aligned}\right.
\end{equation}
where $t \in [0,T)$ is the initial time and $x\in L^2_{\mathcal{F}_t}(\Omega;\mathbb{R}^n)$ is the given initial state. $A(\cdot)$, $B(\cdot)$, $C(\cdot)$, $D(\cdot)$ are given deterministic matrix-valued functions of proper dimensions and $F(\cdot,\cdot), G(\cdot,\cdot)$ are independent of $\omega$. $b(\cdot)$, $\sigma(\cdot)$ are $\mathbb{F}$-progressively measurable processes and $f(\cdot,\cdot)$ is also related to $\omega$. $u(\cdot)$ is the \textit{control process}. We define the admissible control set:
\begin{equation}\label{ac}
\begin{aligned}
\mathcal{U}[t,T]=\bigg\{u:\,&[t,T] \times \Omega \to \mathbb{R}^{m} \,\,\big|\,\,u\,\,\mbox{is}\,\,\mathbb{F}\mbox{-predictable  processes},\,\,\mathbb{E}\int_{t}^{T} |u(s)|^2 ds < \infty \bigg\}.
\end{aligned}
\end{equation}
The control process $u(\cdot)\in\mathcal{U}[t,T]$ is called an \textit{admissible control}.

Then we define the \textit{cost functional}:
\begin{equation}\label{cf}
\begin{aligned}
&J(t,x;u(\cdot))=\mathbb{E}
\Bigg\{\int_t^T \bigg[\big\langle Q(s)X(s),X(s)\big\rangle+2\big\langle S(s)X(s),u(s)\big\rangle+\big\langle R(s)u(s),u(s)\big\rangle\\
&\qquad\qquad +2\big\langle q(s),X(s)\big\rangle+2\big\langle \rho(s),u(s)\big\rangle\bigg] ds+\big\langle HX(T),X(T)\big\rangle+2\big\langle g,X(T)\big\rangle \Bigg\},
\end{aligned}
\end{equation}
where $H$ is a symmetric matrix and $Q(\cdot),S(\cdot)$ and $R(\cdot)$ are deterministic matrix-valued functions of proper dimensions with
$Q(\cdot)^\top=Q(\cdot),\ R(\cdot)^\top=R(\cdot)$. $g$ is an $\mathcal{F}_T$-measurable random variable and $q(\cdot)$ is an $\mathbb{F}$-progressively measurable process and $\rho(\cdot)$ is an $\mathbb{F}$-predictable process.

We adopt the following assumptions.

\textbf{(H1)} The coefficients of the state equation (\ref{state}) satisfy the following:
\begin{equation}\nonumber
\begin{cases}
A(\cdot) \in L^\infty(t,T;\mathbb{R}^{n \times n}),\,\,B(\cdot) \in L^\infty(t,T;\mathbb{R}^{n \times m}),\,\,b(\cdot) \in L^2_\mathbb{F}(t,T;\mathbb{R}^n),\\
C(\cdot) \in L^\infty(t,T;\mathbb{R}^{n \times n}),\,\,D(\cdot) \in L^\infty(t,T;\mathbb{R}^{n \times m}),\,\,\sigma(\cdot) \in L^2_\mathbb{F}(t,T;\mathbb{R}^n),\\
F(\cdot,\cdot) \in F^\infty_{p}(t,T;\mathbb{R}^{n \times n}),\,\,G(\cdot,\cdot) \in  F^\infty_{p}(t,T;\mathbb{R}^{n \times m}),\,\,f(\cdot,\cdot) \in F^2_{\mathbb{F},p}(t,T;\mathbb{R}^n).\\
\end{cases}
\end{equation}

\textbf{(H2)} The weighting coefficients of the cost functionals (\ref{cf}) satisfy the following:
\begin{equation}\nonumber
\begin{cases}
Q(\cdot) \in L^\infty(t,T;\mathbb{S}^n),\,\,S(\cdot) \in L^\infty(t,T;\mathbb{R}^{m \times n}),\,\,R(\cdot) \in L^{\infty}(t,T;\mathbb{S}^{m}),\\
q(\cdot) \in L^2_\mathbb{F}(t,T;\mathbb{R}^n),\,\,\rho(\cdot) \in L^2_{\mathbb{F},p}(t,T;\mathbb{R}^{m}),\,\,g \in L^2_{\mathcal{F}_T}(\Omega;\mathbb{R}^{n}),\,\,H \in \mathbb{S}^n.\\
\end{cases}
\end{equation}

Under (H1), (H2), for any given $x\in L^2_{\mathcal{F}_t}(\Omega;\mathbb{R}^n) $ and $u(\cdot) \in \mathcal{U}[t,T] \equiv L^2_{\mathbb{F},p}(t,T;\mathbb{R}^m)$, the state equation (\ref{state}) admits a unique adapted solution $X(\cdot) \in L^2_{\mathbb{F}}(t,T;\mathbb{R}^n)$ and the cost functional is well-defined. Therefore, the following problem is meaningful.

\textbf{Problem (SLQP)}. For given initial pair $(t,x) \in [0,T]\times L^2_{\mathcal{F}_t}(\Omega;\mathbb{R}^n)$, find a $\bar{u}(\cdot) \in \mathcal{U}[t,T]$ such that
\begin{equation}\label{Problem SLQP}
J(t,x;\bar{u}(\cdot))=\mathop{\inf}\limits_{u(\cdot) \in\, \mathcal{U}[t,T]}J(t,x;u(\cdot)) \equiv V(t,x).
\end{equation}
Any $\bar{u}(\cdot) \in \mathcal{U}[0,T]$ satisfying (\ref{Problem SLQP}) is called \textit{an open-loop optimal control} of Problem (SLQP) for $(t,x)$, the corresponding $\bar{X}(\cdot)\equiv X(\cdot;t,x,\bar{u}(\cdot))$ is called \textit{an open-loop optimal state} and $(\bar{X}(\cdot),\bar{u}(\cdot))$ is called \textit{an open-loop optimal pair}. The map $V(\cdot,\cdot)$ is called \textit{the value function} of Problem (SLQP).

In particular, when $b(\cdot),\,\sigma(\cdot),\,f(\cdot,\cdot),\,q(\cdot),\,\rho(\cdot)$ and $g$ are all zero, we refer to the above problem as the Problem (SLQP)$^0$.

\begin{mydef}\label{def2.1}
For given initial pair $(t,x) \in [0,T]\times L^2_{\mathcal{F}_t}(\Omega;\mathbb{R}^n)$, if there exists a (unique) $\bar{u}(\cdot) \in \mathcal{U}[t,T]$ such that (\ref{Problem SLQP}) holds, then we say that Problem (SLQP) is (uniquely) open-loop solvable for $(t,x)$.
\end{mydef}

Next, take $\Theta(\cdot) \in L^2(t,T;\mathbb{R}^{m \times n}) \equiv \mathcal{Q}[t,T]$ and $v(\cdot) \in \mathcal{U}[t,T]$. For given initial pair $(t,x) \in [0,T)\times L^2_{\mathcal{F}_t}(\Omega;\mathbb{R}^n)$, let us consider the following equation (Some time variables are usually omitted.):
\begin{equation}\label{cl-state}\left\{
\begin{aligned}
dX&=\big[(A+B\Theta)X+Bv+b\big]ds+\big[(C+D\Theta)X+Dv+\sigma\big]dW\\
&\qquad+\int_\mathbf{E}\big[(F(e)+G(e)\Theta)X_{-}+G(e)v+f(e)\big]\tilde{\mathnormal{N}}(deds),\quad s \in [t,T],\\
X(t)&=x,
\end{aligned}\right.
\end{equation}
which admits a unique solution $X(\cdot) \equiv X(\cdot;t,x,\Theta(\cdot),v(\cdot))$, depending on the $\Theta(\cdot)$ and $v(\cdot)$. $(\Theta(\cdot),v(\cdot)) \in \mathcal{Q}[t,T]\times \mathcal{U}[t,T]$ is called a {\it closed-loop strategy} and the above equation (\ref{cl-state}) is called a {\it closed-loop system} of the original state equation (\ref{state}) under $(\Theta(\cdot),v(\cdot))$. We point out that $(\Theta(\cdot),v(\cdot))$ is independent of the initial state $x \in L^2_{\mathcal{F}_t}(\Omega;\mathbb{R}^n)$. With the above solution $X(\cdot)$, we define

\begin{equation}\label{cl-cf}
\begin{aligned}
&J(t,x;\Theta(\cdot)X(\cdot)+v(\cdot))=\mathbb{E} \bigg\{
\int_t^T \bigg[\bigg\langle
\left( \begin{array}{cc} Q & S^\top \\ S & R\end{array} \right)
\left( \begin{array}{c}  X \\ \Theta X+v\end{array} \right),
\left( \begin{array}{c}  X \\ \Theta X+v\end{array}\right)\bigg\rangle\\
&\qquad\qquad\qquad  +2\bigg\langle
\left( \begin{array}{c} q \\ \rho \end{array} \right),
\left( \begin{array}{c} X \\ \Theta X +v \end{array}\right)\bigg\rangle \bigg]ds+\big\langle HX(T),X(T)\big\rangle+2\big\langle g,X(T)\big\rangle \bigg\}\\
&\quad=\mathbb{E}\bigg\{ \int_t^T \bigg[\big\langle (Q+\Theta^\top S+S^\top\Theta+\Theta^\top R\Theta)X,X\big\rangle+2\big\langle (S+R\Theta)X,v\big\rangle+\big\langle Rv,v\big\rangle\\
&\qquad\qquad\qquad+2\big\langle q+\Theta^\top\rho,X\big\rangle+2\big\langle \rho,v\big\rangle\bigg] ds+\big\langle HX(T),X(T)\big\rangle+2\big\langle g,X(T)\big\rangle \bigg\}.
\end{aligned}
\end{equation}

We now introduce the following definition.
\begin{mydef}\label{def2.2}
	A pair $(\bar{\Theta}(\cdot),\bar{v}(\cdot)) \in \mathcal{Q}[t,T] \times \mathcal{U}[t,T]$ is called a \textit{closed-loop optimal strategy} of Problem (SLQP) on $[t,T]$ if
	\begin{equation}
	\begin{aligned}
	&J(t,x;\bar{\Theta}(\cdot)\bar{X}(\cdot)+\bar{v}(\cdot)) \leqslant J(t,x;\Theta(\cdot)X(\cdot)+v(\cdot)),\\
	&\qquad\qquad  \forall x \in L^2_{\mathcal{F}_t}(\Omega;\mathbb{R}^n),\,\, \Theta(\cdot) \in \mathcal{Q}[t,T],\,\, v(\cdot) \in \mathcal{U}[t,T],
	\end{aligned}
	\end{equation}
	where $\bar{X}(\cdot)\equiv X(\cdot;t,x,\bar{\Theta}(\cdot),\bar{v}(\cdot))$ on the left and $X(\cdot)\equiv X(\cdot;t,x,\Theta(\cdot),v(\cdot))$ on the right. In the above case, we say that Problem (SLQP) is {\it closed-loop solvable} on $[t,T]$.
\end{mydef}

We emphasize that the pair $(\bar{\Theta}(\cdot),\bar{v}(\cdot))$ is required to be independent of the initial state $x \in  L^2_{\mathcal{F}_t}(\Omega;\mathbb{R}^n)$. We have the following equivalence theorem.

\begin{mythm}\label{relation}
	Let (H1)-(H2) hold and let $(\bar{\Theta}(\cdot),\bar{v}(\cdot)) \in \mathcal{Q}[t,T] \times \mathcal{U}[t,T]$. Then the following statements are equivalent:\\
	(\romannumeral 1) $(\bar{\Theta}(\cdot),\bar{v}(\cdot))$ is a closed-loop optimal strategy of Problem (SLQP) on $[t,T]$.\\
	(\romannumeral 2) For any given $x \in L^2_{\mathcal{F}_t}(\Omega;\mathbb{R}^n)$ and $v(\cdot) \in \mathcal{U}[t,T]$,
	\begin{equation}\nonumber
	J(t,x;\bar{\Theta}(\cdot)\bar{X}(\cdot)+\bar{v}(\cdot)) \leqslant J(t,x;\bar{\Theta}(\cdot)X(\cdot)+v(\cdot)),
	\end{equation}
	where $X(\cdot)\equiv X(\cdot;t,x,\bar{\Theta}(\cdot),v(\cdot))$ on the right.\\
	(\romannumeral 3) For any given $x \in L^2_{\mathcal{F}_t}(\Omega;\mathbb{R}^n)$ and $u(\cdot) \in \mathcal{U}[t,T]$,
	\begin{equation}\label{open and close relation}
	J(t,x;\bar{\Theta}(\cdot)\bar{X}(\cdot)+\bar{v}(\cdot)) \leqslant J(t,x;u(\cdot)).
	\end{equation}	
\end{mythm}

{\it Proof.} (\romannumeral 1) $\Rightarrow$ (\romannumeral 2). From the definition of closed-loop optimal strategy, it can be proved.\\
(\romannumeral 2) $\Rightarrow$ (\romannumeral 3). For given $x \in L^2_{\mathcal{F}_t}(\Omega;\mathbb{R}^n)$ and $u(\cdot) \in \mathcal{U}[t,T]$, $X(\cdot)$ is the adapted solution to the following SDEP:
\begin{equation*}\left\{
\begin{aligned}
dX(s)&=\big[A(s)X(s)+B(s)u(s)+b(s)\big]ds+\big[C(s)X(s)+D(s)u(s)+\sigma(s)\big]dW(s)\\
&\qquad+\int_\mathbf{E}\big[F(s,e)X(s-)+G(s,e)u(s)+f(s,e)\big]\tilde{\mathnormal{N}}(deds),\quad s \in [t,T],\\
X(t)&=x,
\end{aligned}\right.
\end{equation*}
Taking $v(\cdot)=u(\cdot)-\bar{\Theta}(\cdot)X(\cdot)$, it is easy to see that $v(\cdot) \in \mathcal{U}[t,T]$. Thus
\begin{equation*}
u(\cdot)=\bar{\Theta}(\cdot)X(\cdot)+v(\cdot),
\end{equation*}
then
\begin{equation*}
X(\cdot)=X(\cdot;t,x,u(\cdot))=X(\cdot;t,x,\bar{\Theta}(\cdot),v(\cdot)) \equiv \tilde{X}(\cdot).
\end{equation*}
Therefore, by (\romannumeral 2), one has
\begin{equation*}
J(t,x;\bar{\Theta}(\cdot)\bar{X}(\cdot)+\bar{v}(\cdot)) \leqslant J(t,x;\bar{\Theta}(\cdot)\tilde{X}(\cdot)+v(\cdot))=J(t,x;u(\cdot)),
\end{equation*}
proving (\romannumeral 3).\\
(\romannumeral 3) $\Rightarrow$ (\romannumeral 1). For a given $x \in L^2_{\mathcal{F}_t}(\Omega;\mathbb{R}^n)$,\,\,$u(\cdot) \in \mathcal{U}[t,T]$ and $(\Theta(\cdot),v(\cdot)) \in \mathcal{Q}[t,T] \times \mathcal{U}[t,T]$, let $X(\cdot) \equiv X(\cdot;t,x,\Theta(\cdot),v(\cdot))$. Taking
\begin{equation*}
u(\cdot)=\Theta(\cdot)X(\cdot)+v(\cdot),
\end{equation*}
from the existence and uniqueness of the adapt solution of SDEP, we know that
\begin{equation*}
X(\cdot;t,x,\Theta(\cdot),v(\cdot))=X(\cdot;t,x,u(\cdot)).
\end{equation*}
Therefore, by (\romannumeral 3), we have
\begin{equation*}
J(t,x;\bar{\Theta}(\cdot)\bar{X}(\cdot)+\bar{v}(\cdot)) \leqslant J(t,x;u(\cdot))= J(t,x;\Theta(\cdot)X(\cdot)+v(\cdot)).
\end{equation*}
This completes the proof. $\qquad\Box$

From the above result, we see that if $(\bar{\Theta}(\cdot),\bar{v}(\cdot))$ is a closed-loop optimal strategy of Problem (SLQP) on $[t,T]$, then for any given initial state $x \in L^2_{\mathcal{F}_t}(\Omega;\mathbb{R}^n)$ and $u(\cdot) \in \mathcal{U}[t,T]$, with the solution $\bar{X}(\cdot)=\bar{X}(\cdot;t,x,\bar{\Theta}(\cdot),\bar{v}(\cdot))$, (\ref{open and close relation}) implies that the {\it outcome}
\begin{equation}\nonumber
\bar{u}(\cdot)=\bar{\Theta}(\cdot)\bar{X}(\cdot)+\bar{v}(\cdot) \in \mathcal{U}[t,T]
\end{equation}
is an open-loop optimal control of Problem (SLQP) for $x$. Therefore, for Problem (SLQP), the existence of closed-loop optimal strategies on $[t,T]$ implies the existence of open-loop optimal controls for any $x \in L^2_{\mathcal{F}_t}(\Omega;\mathbb{R}^n)$.

The following result is concerned with open-loop solvability of Problem (SLQP) for given initial state.

% SLQP open-loop optimal solvability
\begin{mypro}\label{opsn}
Let (H1)-(H2) hold. For given initial pair $(t,x) \in [0,T]\times L^2_{\mathcal{F}_t}(\Omega;\mathbb{R}^n)$, a control $\bar{u}(\cdot)$ is an open-loop optimal control of Problem (SLQP) if and only if the following hold:\\
(\romannumeral 1) The stationarity condition holds:
\begin{equation}\label{sc}
B^\top\bar{Y}+D^\top\bar{Z}+\int_\mathbf{E}G(e)^\top\bar{K}(e)\pi(de)+S\bar{X}+R\bar{u}+\rho=0,\quad a.e.,\ \mathbb{P}\mbox{-}a.s.,
\end{equation}
where $(\bar{X}(\cdot),\bar{Y}(\cdot),\bar{Z}(\cdot),\bar{K}(\cdot,\cdot)) \in L^2_{\mathbb{F}}(t,T;\mathbb{R}^n) \times L^2_{\mathbb{F}}(t,T;\mathbb{R}^n) \times L^2_{\mathbb{F},p}(t,T;\mathbb{R}^n)\times F^2_{\mathbb{F},p}(t,T;\mathbb{R}^n) $ is the adapted solution to the following FBSDEP:
\begin{equation}\left\{
\begin{aligned}
d\bar{X}&=\big[A\bar{X}+B\bar{u}+b\big]ds+\big[C\bar{X}+D\bar{u}+\sigma\big]dW\\
&\qquad+\int_\mathbf{E}\big[F(e)\bar{X}_{-}+G(e)\bar{u}+f(e)\big]\tilde{\mathnormal{N}}(deds),\\
d\bar{Y}&=-\big[A^\top\bar{Y}+C^\top\bar{Z}+\int_\mathbf{E}F^\top(e)\bar{K}(e)\pi(de)+Q\bar{X}+S^\top\bar{u}+q\big]ds\\
&\qquad+\bar{Z}dW+\int_\mathbf{E}\bar{K}(e)\tilde{\mathnormal{N}}(deds),\quad s \in [t,T],\\
\bar{X}(t)&=x,\quad \bar{Y}(T)=H\bar{X}(T)+g.
\end{aligned}\right.
\end{equation}
(\romannumeral 2) The convexity condition holds: For any $u(\cdot) \in \mathcal{U}[t,T]$,
\begin{equation}\label{cc}
\mathbb{E}\bigg[\big\langle Hx_0(T),x_0(T)\big\rangle+\int_t^T \Big(\big\langle Qx_0,x_0 \big\rangle +2\big\langle Sx_0,u \big\rangle +\big\langle Ru,u \big\rangle \Big)ds\bigg] \geqslant 0,
\end{equation}
where $x_0(\cdot) \in L^2_{\mathbb{F}}(t,T;\mathbb{R}^n)$ is the adapted solution to the following SDEP:
\begin{equation}\left\{
\begin{aligned}
dx_0(s)&=\big[A(s)x_0(s)+B(s)u(s)\big]ds+\big[C(s)x_0(s)+D(s)u(s)\big]dW(s)\\
&\qquad+\int_\mathbf{E}\big[F(s,e)x_0(s-)+G(s,e)u(s)\big]\tilde{\mathnormal{N}}(deds),\quad s \in [t,T],\\
x_0(t)&=0.
\end{aligned}\right.
\end{equation}
\end{mypro}

\textit{Proof.} For any $u(\cdot) \in \mathcal{U}[t,T]$ and $\epsilon \in \mathbb{R}$, let $u^\epsilon(\cdot)=\bar{u}(\cdot)+\epsilon u(\cdot)$, thus $X^{\epsilon}(\cdot)\equiv X(\cdot;t,x,\bar{u}(\cdot)+\epsilon u(\cdot))$ is the corresponding state which satisfies
\begin{equation*}
\left\{
\begin{aligned}
dX^{\epsilon}&=\big[AX^{\epsilon}+B(\bar{u}+\epsilon u)+b\big]ds+\big[CX^{\epsilon}+D(\bar{u}+\epsilon u)+\sigma\big]dW\\
&\qquad+\int_\mathbf{E}\big[F(e)X^{\epsilon}_{-}+G(e)(\bar{u}+\epsilon u)+f(e)\big]\tilde{\mathnormal{N}}(deds),\quad s \in [t,T],\\
X^{\epsilon}(t)&=x.
\end{aligned}\right.
\end{equation*}
Then $\frac{X^{\epsilon}(\cdot)-X(\cdot)}{\epsilon}$ satisfies the following SDEP:
\begin{equation*}
\left\{
\begin{aligned}
d\bigg(\frac{X^{\epsilon}-X}{\epsilon}\bigg)&=\bigg[A\bigg(\frac{X^{\epsilon}-X}{\epsilon}\bigg)+Bu\bigg]ds+\bigg[C\bigg(\frac{X^{\epsilon}-X}{\epsilon}\bigg)+Du\bigg]dW\\
&\qquad+\int_\mathbf{E}\bigg[F(e)\bigg(\frac{X^{\epsilon}_{-}-X_{-}}{\epsilon}\bigg)+G(e)u\bigg]\tilde{\mathnormal{N}}(deds),\quad s \in [t,T],\\
\frac{X^{\epsilon}-X}{\epsilon}(t)&=0.
\end{aligned}\right.
\end{equation*}
From the existence and uniqueness of the solution to SDEP, we know that $x_0 \equiv \frac{X^{\epsilon}-X}{\epsilon}$. Then
\begin{equation*}
\begin{aligned}
&J(t,x;\bar{u}(\cdot)+\epsilon u(\cdot))-J(t,x;\bar{u}(\cdot))\\
&=2\epsilon\mathbb{E}\biggl\{\big\langle H\bar{X}(T)+g,x_0(T)\big\rangle
+\int_t^T\Big[ \big\langle Qx_0,\bar{X}\big\rangle +\big\langle S\bar{X},u\big\rangle \\
&\qquad\qquad +\big\langle Sx_0,\bar{u} \big\rangle+\big\langle R\bar{u},u \big\rangle +\big\langle q,x_0\big\rangle+\big\langle \rho,u \big\rangle \Big]ds\biggr\}\\
&\quad\,+\epsilon^2\mathbb{E}\biggl\{\big\langle Hx_0(T),x_0(T)\rangle+\int_t^T\Big[ \big\langle Qx_0,x_0 \big\rangle+2\big\langle Sx_0,u\big\rangle
+\big\langle Ru,u\big\rangle \Big]ds \biggr\}.
\end{aligned}
\end{equation*}
Applying It\^o's formula to $\big\langle \bar{Y}(\cdot),x_0(\cdot) \big\rangle$, we get
\begin{equation*}
\begin{aligned}
&J(t,x;\bar{u}(\cdot)+\epsilon u(\cdot))-J(t,x;\bar{u}(\cdot))\\
&=2\epsilon\mathbb{E}\int_t^T \big\langle B^\top\bar{Y}+D^\top\bar{Z}+\int_\mathbf{E}G(e)^\top\bar{K}(e)\pi(de)+S\bar{X}+R\bar{u}+\rho,u \big\rangle ds\\
&\quad\,+\epsilon^2\mathbb{E}\biggl\{\big\langle Hx_0(T),x_0(T)\rangle+\int_t^T\Big[ \big\langle Qx_0,x_0 \big\rangle+2\big\langle Sx_0,u\big\rangle
+\big\langle Ru,u\big\rangle \Big]ds \biggr\}.
\end{aligned}
\end{equation*}
Therefore, $(\bar{X}(\cdot),\bar{u}(\cdot))$ is an open-loop optimal pair of Problem (SLQP) if and only if (\ref{sc}) and (\ref{cc}) hold. The proof is complete. $\qquad\Box$

\vspace{2mm}

On the other hand, from the second part of Theorem \ref{relation}, we can see that $(\bar{\Theta}(\cdot),\bar{v}(\cdot))$ being a closed-loop optimal strategy of Problem (SLQP) is equivalent to $\bar{v}(\cdot)$ being an open-loop optimal control of the SLQP optimal control problem (\ref{cl-state})-(\ref{cl-cf}) with $\Theta(\cdot)=\bar{\Theta}(\cdot)$, which we denote by {\bf Problem (SLQP)$_{\bar{\Theta}}$}. In particular, when $b(\cdot),\,\sigma(\cdot),\,f(\cdot,\cdot),\,q(\cdot),\,\rho(\cdot)$ and $g$ are all zero, we refer to it as Problem (SLQP)$_{\bar{\Theta}}^0$.
Similar to Proposition \ref{opsn}, we can give the following result.

\begin{mypro}\label{SLQP_bar(Theta)}
Let (H1)-(H2) hold. For any given $x \in L^2_{\mathcal{F}_t}(\Omega;\mathbb{R}^n)$, $\bar{v}(\cdot)$ is an open-loop optimal control of Problem (SLQP)$_{\bar{\Theta}}$ if and only if the following stationarity condition holds:
\begin{equation}\label{cl-sc}
B^\top\bar{Y}+D^\top\bar{Z}+\int_\mathbf{E}G(e)^\top\bar{K}(e)\pi(de)+(S+R\bar{\Theta})\bar{X}+R\bar{v}+\rho=0,\quad a.e.,\ \mathbb{P}\mbox{-}a.s.,
\end{equation}
where $(\bar{X}(\cdot),\bar{Y}(\cdot),\bar{Z}(\cdot),\bar{K}(\cdot,\cdot)) \in L^2_{\mathbb{F}}(t,T;\mathbb{R}^n) \times L^2_{\mathbb{F}}(t,T;\mathbb{R}^n) \times L^2_{\mathbb{F},p}(t,T;\mathbb{R}^n)\times F^2_{\mathbb{F},p}(t,T;\mathbb{R}^n) $ is the adapted solution to the following FBSDEP:
\begin{equation}\label{cl-os-1}\left\{
\begin{aligned}
d\bar{X}&=\big[(A+B\bar{\Theta})\bar{X}+B\bar{v}+b\big]ds+\big[(C+D\bar{\Theta})\bar{X}+D\bar{v}+\sigma\big]dW\\
&\qquad+\int_\mathbf{E}\big[(F(e)+G(e)\bar{\Theta})\bar{X}_{-}+G(e)\bar{v}+f(e)\big]\tilde{\mathnormal{N}}(deds),\\
d\bar{Y}&=-\bigg[(A+B\bar{\Theta})^\top\bar{Y}+(C+D\bar{\Theta})^\top\bar{Z}+\int_\mathbf{E}(F(e)+G(e)\bar{\Theta})^\top\bar{K}(e)\pi(de)\\
&\qquad\quad+(Q+S^\top\bar{\Theta}+\bar{\Theta}^\top S+\bar{\Theta}^\top R\bar{\Theta})\bar{X}+(S+R\bar{\Theta})^\top\bar{v}+q+\bar{\Theta}^\top\rho\bigg]ds\\
&\qquad+\bar{Z}dW+\int_\mathbf{E}\bar{K}(e)\tilde{\mathnormal{N}}(deds),\quad s \in [t,T],\\
\bar{X}(t)&=x,\quad \bar{Y}(T)=H\bar{X}(T)+g,
\end{aligned}\right.
\end{equation}
and the following convexity condition holds: For any $v(\cdot) \in \mathcal{U}[t,T]$,
\begin{equation*}
\begin{aligned}
&\mathbb{E}\bigg\{ \int_t^T \Big[\big\langle (Q+\bar{\Theta}^\top S+S^\top\bar{\Theta}+\bar{\Theta}^\top R\bar{\Theta})X,X\big\rangle+2\big\langle (S+R\bar{\Theta})X,v\big\rangle\\
&\qquad\qquad\qquad  +\big\langle Rv,v\big\rangle\Big] ds+\big\langle HX(T),X(T)\big\rangle \bigg\}\ge 0,
\end{aligned}
\end{equation*}
where $X(\cdot)\in L^2_{\mathbb{F}}(t,T;\mathbb{R}^n)$ is the adapted solution to the following SDEP:
\begin{equation*}\left\{
\begin{aligned}
dX&=\big[(A+B\bar{\Theta})X+Bv\big]ds+\big[(C+D\bar{\Theta})X+Dv\big]dW\\
&\qquad+\int_\mathbf{E}\big[(F(e)+G(e)\bar{\Theta})X_{-}+G(e)v\big]\tilde{\mathnormal{N}}(deds),\quad s \in [t,T],\\
X(t)&=0.
\end{aligned}\right.
\end{equation*}
\end{mypro}

\section{Main results}

In this section, we will study the necessary and sufficient conditions for Problem (SLQP) to be closed-loop solvable.

Making use of (\ref{cl-sc}), we may rewrite the BSDEP in (\ref{cl-os-1}) and obtain
%as follows:
%\begin{equation*}
%\begin{aligned}
%d\bar{Y}%&=-\bigg[(A+B\bar{\Theta})^\top\bar{Y}+(C+D\bar{\Theta})^\top\bar{Z}+\int_\mathbf{E}(F(e)+G(e)\bar{\Theta})^\top\bar{K}(e)\pi(de)\\
%&\qquad\quad+(Q+S^\top\bar{\Theta}+\bar{\Theta}^\top S+\bar{\Theta}^\top R\bar{\Theta})\bar{X}+(S+R\bar{\Theta})^\top\bar{v}+q+\bar{\Theta}^\top\rho\bigg]ds\\
%&\quad\quad+\bar{Z}dW+\int_\mathbf{E}\bar{K}(e)\tilde{\mathnormal{N}}(deds)\\
%&=-\bigg[A^\top\bar{Y}+C^\top\bar{Z}+\int_\mathbf{E}F(e)^\top\bar{K}(e)\pi(de)+(Q+S^\top\bar{\Theta})\bar{X}+S^\top\bar{v}+q\\
%&\qquad\quad+\bar{\Theta}^\top\big[ B^\top\bar{Y}+D^\top\bar{Z}+\int_\mathbf{E}G(e)^\top\bar{K}(e)\pi(de)+(S+R\bar{\Theta})\bar{X}+R\bar{v}+\rho \big]\bigg]ds\\
%&\qquad+\bar{Z}dW+\int_\mathbf{E}\bar{K}(e)\tilde{\mathnormal{N}}(deds)\\
%&=-\bigg[A^\top\bar{Y}+C^\top\bar{Z}+\int_\mathbf{E}F(e)^\top\bar{K}(e)\pi(de)+(Q+S^\top\bar{\Theta})\bar{X}+S^\top\bar{v}+q\bigg]ds\\
%&\qquad\quad+\bar{Z}dW+\int_\mathbf{E}\bar{K}(e)\tilde{\mathnormal{N}}(deds).\\
%\end{aligned}
%\end{equation*}
%Thus, we obtain
\begin{equation}\label{cl-os}\left\{
\begin{aligned}
d\bar{X}&=\big[(A+B\bar{\Theta})\bar{X}+B\bar{v}+b\big]ds+\big[(C+D\bar{\Theta})\bar{X}+D\bar{v}+\sigma\big]dW\\
&\qquad\quad+\int_\mathbf{E}\big[(F(e)+G(e)\bar{\Theta})\bar{X}_{-}+G(e)\bar{v}+f(e)\big]\tilde{\mathnormal{N}}(deds),\\
d\bar{Y}&=-\bigg[A^\top\bar{Y}+C^\top\bar{Z}+\int_\mathbf{E}F(e)^\top\bar{K}(e)\pi(de)+(Q+S^\top\bar{\Theta})\bar{X}+S^\top\bar{v}+q\bigg]ds\\
&\qquad+\bar{Z}dW+\int_\mathbf{E}\bar{K}(e)\tilde{\mathnormal{N}}(deds),\quad s \in [t,T],\\
\bar{X}(t)&=x,\quad \bar{Y}(T)=H\bar{X}(T)+g.
\end{aligned}\right.
\end{equation}

We first have the following result.

\begin{mythm}\label{problem transform}
Let (H1)-(H2) hold. If $(\bar{\Theta}(\cdot),\bar{v}(\cdot)) \in \mathcal{Q}[t,T] \times \mathcal{U}[t,T]$ is an optimal closed-loop strategy of Problem (SLQP) on $[t,T]$, then $(\bar{\Theta}(\cdot),0)$ is an optimal closed-loop strategy of Poblem (SLQP)$^0$ on $[t,T]$.	
\end{mythm}
	
\textit{Proof.} By the Proposition (\ref{SLQP_bar(Theta)}), we can see that $(\bar{\Theta}(\cdot),\bar{v}(\cdot))$ is an optimal closed-loop strategy of Problem (SLQP) on $[t,T]$ if and only if for any $x \in L^2_{\mathcal{F}_t}(\Omega;\mathbb{R}^n)$, the stationarity condition holds:
\begin{equation}\label{Th3.1-sc}
B^\top\bar{Y}+D^\top\bar{Z}+\int_\mathbf{E}G(e)^\top\bar{K}(e)\pi(de)+(S+R\bar{\Theta})\bar{X}+R\bar{v}+\rho=0,\quad a.e.,\ \mathbb{P}\mbox{-}a.s.,
\end{equation}
where $(\bar{X}(\cdot),\bar{Y}(\cdot),\bar{Z}(\cdot),\bar{K}(\cdot,\cdot))$ is the adapted solution to the following FBSDEP:
\begin{equation}\label{Th3.1-os}\left\{
\begin{aligned}
d\bar{X}&=\big[(A+B\bar{\Theta})\bar{X}+B\bar{v}+b\big]ds+\big[(C+D\bar{\Theta})\bar{X}+D\bar{v}+\sigma\big]dW\\
&\qquad+\int_\mathbf{E}\big[(F(e)+G(e)\bar{\Theta})\bar{X}_{-}+G(e)\bar{v}+f(e)\big]\tilde{\mathnormal{N}}(deds),\\
d\bar{Y}&=-\bigg[(A+B\bar{\Theta})^\top\bar{Y}+(C+D\bar{\Theta})^\top\bar{Z}+\int_\mathbf{E}(F(e)+G(e)\bar{\Theta})^\top\bar{K}(e)\pi(de)\\
&\qquad\quad+(Q+S^\top\bar{\Theta}+\bar{\Theta}^\top S+\bar{\Theta}^\top R\bar{\Theta})\bar{X}+(S+R\bar{\Theta})^\top\bar{v}+q+\bar{\Theta}^\top\rho\bigg]ds\\
&\quad\quad+\bar{Z}dW+\int_\mathbf{E}\bar{K}(e)\tilde{\mathnormal{N}}(deds),\quad s \in [t,T],\\
\bar{X}(t)&=x,\quad \bar{Y}(T)=H\bar{X}(T)+g,
\end{aligned}\right.
\end{equation}
and the following convexity condition holds: For any $v(\cdot) \in \mathcal{U}[t,T]$,
\begin{equation*}
\begin{aligned}
&\mathbb{E}\bigg\{ \int_t^T \bigg[\big\langle (Q+\bar{\Theta}^\top S+S^\top\bar{\Theta}+\bar{\Theta}^\top R\bar{\Theta})X,X\big\rangle+2\big\langle (S+R\bar{\Theta})X,v\big\rangle\\
&\qquad\qquad\qquad  +\big\langle Rv,v\big\rangle\bigg] ds+\big\langle HX(T),X(T)\big\rangle \bigg\}\ge 0,
\end{aligned}
\end{equation*}
where $X(\cdot)$ is the adapted solution to the following SDEP:
\begin{equation*}\left\{
\begin{aligned}
dX&=\big[(A+B\bar{\Theta})X+Bv\big]ds+\big[(C+D\bar{\Theta})X+Dv\big]dW\\
&\qquad+\int_\mathbf{E}\big[(F(e)+G(e)\bar{\Theta})X_{-}+G(e)v\big]\tilde{\mathnormal{N}}(deds),\quad s \in [t,T],\\
X(t)&=0,
\end{aligned}\right.
\end{equation*}

Since the FBSDEP (\ref{Th3.1-os}) admits a solution for each $x \in L^2_{\mathcal{F}_t}(\Omega;\mathbb{R}^n)$ and $(\bar{\Theta}(\cdot),\bar{v}(\cdot))$ is independent of $x$, by subtracting solutions corresponding $x$ and $0$, the latter from the former, we see that for any $x \in L^2_{\mathcal{F}_t}(\Omega;\mathbb{R}^n)$, the following FBSDEP:
\begin{equation}\label{Th3.1-os-2}\left\{
\begin{aligned}
dX&=\big[(A+B\bar{\Theta})X\big]ds+\big[(C+D\bar{\Theta})X\big]dW+\int_\mathbf{E}\big[(F(e)+G(e)\bar{\Theta})X_{-}\big]\tilde{\mathnormal{N}}(deds),\\
dY&=-\bigg[(A+B\bar{\Theta})^\top Y+(C+D\bar{\Theta})^\top Z+\int_\mathbf{E}(F(e)+G(e)\bar{\Theta})^\top K(e)\pi(de)\\
&\qquad+(Q+S^\top\bar{\Theta}+\bar{\Theta}^\top S+\bar{\Theta}^\top R\bar{\Theta})X\bigg]ds+ZdW+\int_\mathbf{E}K(e)\tilde{\mathnormal{N}}(deds),\quad s \in [t,T],\\
X(t)&=x,\quad Y(T)=HX(T),
\end{aligned}\right.
\end{equation}
ant in this time $(X(\cdot),Y(\cdot),Z(\cdot),K(\cdot,\cdot))$ satisfies
\begin{equation*}
B^\top Y+D^\top Z+\int_\mathbf{E}G(e)^\top K(e)\pi(de)+(S+R\bar{\Theta})X=0,\quad a.e.,\ \mathbb{P}\mbox{-}a.s..
\end{equation*}
It follows, again from Theorem \ref{relation} and Proposition \ref{SLQP_bar(Theta)}, that $(\bar{\Theta}(\cdot),0)$ is an optimal closed-loop strategy of Problem (SLQP)$^0$ on $[t,T]$.
The proof is complete. $\qquad\Box$  	

To summarize the relationship between Problem (SLQP), Problem (SLQP)$^0$, Problem (SLQP)$_{\bar{\Theta}}$ and Problem (SLQP)$^0_{\bar{\Theta}}$, we plot the following diagram:

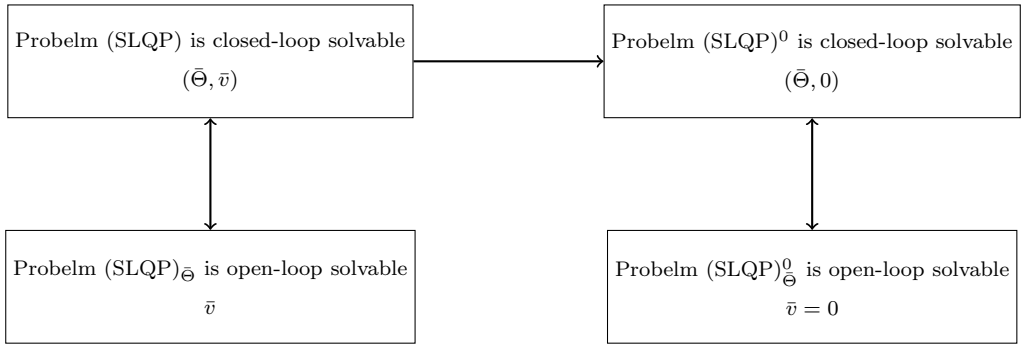
\vspace{2mm}

\begin{figure}[h]
\scriptsize
\tikzstyle{format} = [rectangle, minimum width = 2cm, minimum height=1.5cm, align=center, draw = black]
%\tikzstyle{format}=[rectangle,draw,thick,fill=white]
\tikzstyle{test}=[diamond,aspect=10,draw,thin]
\tikzstyle{point}=[coordinate,on grid,]
\begin{tikzpicture}
\node[format] (SLQP)
{Probelm (SLQP) is closed-loop solvable \\[1ex] $(\bar{\Theta},\bar{v})$};
\node[format,below of=SLQP,node distance=3cm](SLQPT)
{Probelm (SLQP)$_{\bar{\Theta}}$ is open-loop solvable \\[1ex] $\bar{v}$};
\node[format,right of=SLQP,node distance=8cm](SLQP0)
{Probelm (SLQP)$^0$ is closed-loop solvable \\[1ex] $(\bar{\Theta},0)$};
\node[format,below of=SLQP0,node distance=3cm](SLQPT0)
{Probelm (SLQP)$^0_{\bar{\Theta}}$ is open-loop solvable \\[1ex] $\bar{v}=0$};
\draw[->,thick](SLQP)--(SLQP0);
\draw[<->,thick](SLQP)--(SLQPT);
\draw[<->,thick](SLQP0)--(SLQPT0);
\end{tikzpicture}
\end{figure}

It is clear that when we want to study the necessary conditions for the closed-loop solvability of Problem (SLQP), we can transform the original problem into the open-loop solvability of Problem (SLQP)$^0_{\bar{\Theta}}$ where the open-loop optimal control is $\bar{v}(\cdot)\equiv 0$. Thus from Proposition \ref{opsn}, we can know that the optimal system of Problem (SLQP)$^0_{\bar{\Theta}}$ is the following FBSDEP:
\begin{equation}\label{SLQT0_os}\left\{
\begin{aligned}
d\bar{X}&=\big[(A+B\bar{\Theta})\bar{X}\big]ds+\big[(C+D\bar{\Theta})\bar{X}\big]dW\\
&\qquad+\int_\mathbf{E}\big[(F(e)+G(e)\bar{\Theta})\bar{X}_{-}\big]\tilde{\mathnormal{N}}(deds),\\
d\bar{Y}&=-\bigg[A^\top\bar{Y}+C^\top\bar{Z}+\int_\mathbf{E}F(e)^\top\bar{K}(e)\pi(de)+(Q+S^\top\bar{\Theta})\bar{X}\bigg]ds\\
&\qquad+\bar{Z}dW+\int_\mathbf{E}\bar{K}(e)\tilde{\mathnormal{N}}(deds),\quad s \in [t,T],\\
\bar{X}(t)&=x,\quad \bar{Y}(T)=H\bar{X}(T),
\end{aligned}\right.
\end{equation}
ant in this time $(\bar{X}(\cdot),\bar{Y}(\cdot),\bar{Z}(\cdot),\bar{K}(\cdot,\cdot))$ satisfies
\begin{equation}\label{SLQT0_sc}
B^\top \bar{Y}+D^\top \bar{Z}+\int_\mathbf{E}G(e)^\top \bar{K}(e)\pi(de)+(S+R\bar{\Theta})\bar{X}=0,\quad a.e.,\ \mathbb{P}\mbox{-}a.s..
\end{equation}

In the light of $\bar{Y}(T)=H\bar{X}(T)$, we assume that
\begin{equation}\label{y=px}
\bar{Y}(\cdot)=P(\cdot)\bar{X}(\cdot),
\end{equation}
where $P(\cdot): [t,T] \rightarrow \mathbb{S}^n$ is a matrix-valued differential function satisfying $P(T)=H$. Applying It\^o-Wentzell's formula to (\ref{y=px}), we have
\begin{equation}\label{y=px_1}
\begin{aligned}
d\bar{Y}&=\dot{P}\bar{X}ds+Pd\bar{X}\\
&=\big[\dot{P}\bar{X}+P(A+B\bar{\Theta})\bar{X}\big]ds+P(C+D\bar{\Theta})\bar{X}dW+\int_\mathbf{E} P(F(e)+G(e)\bar{\Theta})\bar{X}_{-}\tilde{\mathnormal{N}}(deds),
\end{aligned}
\end{equation}

Comparing the diffusion coefficient of the above equation and of the second equation in (\ref{SLQT0_os}), we can see
\begin{equation}\label{zk}\left\{
\begin{aligned}
&\bar{Z}=P(C+D\bar{\Theta})\bar{X},\\
&\bar{K}(e)= P(F(e)+G(e)\bar{\Theta})\bar{X}.
\end{aligned}\right.
\end{equation}
Plugging (\ref{y=px}) and (\ref{zk}) into (\ref{SLQT0_sc}), we obtain
\begin{equation}\label{sc_1}
\bigg[B^\top P+D^\top P(C+D\bar{\Theta})+\int_\mathbf{E}G(e)^\top P(F(e)+G(e)\bar{\Theta})\pi(de)+(S+R\bar{\Theta})\bigg]\bar{X}=0
\end{equation}	

Now comparing the drift coefficient of (\ref{y=px_1}) and of the second equation in (\ref{SLQT0_os}), noting that (\ref{y=px}), (\ref{zk}) and (\ref{sc_1}), we have
\begin{equation*}
\begin{aligned}
0&=\bigg[\dot{P}+P(A+B\bar{\Theta})+A^\top P+C^\top P(C+D\bar{\Theta})\\
&\qquad+\int_\mathbf{E}F(e)^\top P(F(e)+G(e)\bar{\Theta})\pi(de)+Q+S^\top\bar{\Theta}\bigg]\bar{X}\\
&=\bigg[\dot{P}+P(A+B\bar{\Theta})+A^\top P+C^\top P(C+D\bar{\Theta})+Q+S^\top\bar{\Theta}\\
&\qquad +\int_\mathbf{E}F(e)^\top P(F(e)+G(e)\bar{\Theta})\pi(de)+\bar{\Theta}^\top\bigg(B^\top P+D^\top P(C+D\bar{\Theta})+S\\
&\qquad +R\bar{\Theta}+\int_\mathbf{E}G(e)^\top P(F(e)+G(e)\bar{\Theta})\pi(de)\bigg)\bigg]\bar{X}\\
&=\bigg[\dot{P}+P(A+B\bar{\Theta})+(A+B\bar{\Theta})^\top P+(C+D\bar{\Theta})^\top P(C+D\bar{\Theta})+Q+S^\top\bar{\Theta}\\
&\qquad +\bar{\Theta}^\top S+\bar{\Theta}^\top R\bar{\Theta}+\int_\mathbf{E}(F(e)+G(e)\bar{\Theta})^\top P(F(e)+G(e)\bar{\Theta})\pi(de)\bigg]\bar{X}.
\end{aligned}
\end{equation*}

Thus, we let $P(\cdot)$ satisfy the following Lyapunov integral-differential equation:
\begin{equation}\label{Lyapunov}\left\{
\begin{aligned}
&0=\dot{P}+P(A+B\bar{\Theta})+(A+B\bar{\Theta})^\top P+(C+D\bar{\Theta})^\top P(C+D\bar{\Theta})+Q+S^\top\bar{\Theta}+\bar{\Theta}^\top S\\
&\qquad +\bar{\Theta}^\top R\bar{\Theta}+\int_\mathbf{E}(F(e)+G(e)\bar{\Theta})^\top P(F(e)+G(e)\bar{\Theta})\pi(de),\quad s \in [t,T],\\
&P(T)=0.
\end{aligned}\right.
\end{equation}	
	
\begin{mypro}\label{Condition of regularity}
Let $P(\cdot)$ be the solution to (\ref{Lyapunov}). Then for any $s\in[t,T]$, we have
\begin{equation}\label{regularity}\left\{
\begin{aligned}
&R+D^\top PD+\int_\mathbf{E}G(e)^\top PG(e)\pi(de) \geqslant 0,\\
&0=\bigg[R+D^\top PD+\int_\mathbf{E}G(e)^\top PG(e)\pi(de) \bigg]\bar{\Theta}\\
&\quad\quad +B^\top P+D^\top PC+S+\int_\mathbf{E}G(e)^\top PF(e)\pi(de).
\end{aligned}\right.
\end{equation}
\end{mypro}	
	
\textit{Proof}. Let us consider Problem (SLQP)$^0_{\bar{\Theta}}$. For any $v(\cdot) \in \mathcal{U}[t,T]$, the state equation and the cost functional are:
\begin{equation}\label{slqpt0-state}\left\{
\begin{aligned}
dX&=\big[(A+B\bar{\Theta})X+Bv\big]ds+\big[(C+D\bar{\Theta})X+Dv\big]dW\\
&\qquad+\int_\mathbf{E}\big[(F(e)+G(e)\bar{\Theta})X_{-}+G(e)v\big]\tilde{\mathnormal{N}}(deds),\quad s \in [t,T],\\
X(t)&=x,
\end{aligned}\right.
\end{equation}
and
\begin{equation}\label{slqpt0-cf}
\begin{aligned}
\tilde{J}(t,x;v(\cdot))&=\mathbb{E}
\bigg\{\int_t^T \bigg[\big\langle \big(Q+\bar{\Theta}^\top S+S\bar{\Theta}+\bar{\Theta}^\top R\bar{\Theta}\big)X,X\big\rangle+2\big\langle \big(S+R\bar{\Theta}\big)X,v\big\rangle\\
& \qquad\qquad \qquad  +\big\langle Rv,v\big\rangle\bigg] ds+\big\langle HX(T),X(T)\big\rangle \bigg\}.
\end{aligned}
\end{equation}
Applying It\^o-Wentzell's formula to $\big\langle PX(\cdot),X(\cdot) \big\rangle$, we get
\begin{equation*}
\begin{aligned}
d\big\langle PX,X\rangle&=\Big\langle -P(A+B\bar{\Theta})X-(A+B\bar{\Theta})^\top PX-(C+D\bar{\Theta})^\top P(C+D\bar{\Theta})X\\
&\qquad-\int_\mathbf{E}\big[(F(e)+G(e)\bar{\Theta})^\top P(F(e)+G(e)\bar{\Theta})X\big]\pi(de)\\
&\qquad-(Q+\bar{\Theta}^\top S+S^\top\bar{\Theta}+\bar{\Theta}^\top R\bar{\Theta})X,X\Big\rangle ds+2\big\langle P\big[(A+B\bar{\Theta})X+Bv\big],X\big\rangle ds\\
&\quad+\big\langle P\big[(C+D\bar{\Theta})X+Dv\big],(C+D\bar{\Theta})X+Dv\big\rangle ds\\
&\quad+\int_\mathbf{E}\big\langle P\big[(F(e)+G(e)\bar{\Theta})X+G(e)v\big],(F(e)+G(e)\bar{\Theta})X+G(e)v\big\rangle \pi(de)ds\\
&\quad+[\cdots]dW+\int_\mathbf{E}[\cdots]\tilde{\mathnormal{N}}(deds)\\
&=\bigg\{-\big\langle\big[Q+\bar{\Theta}^\top S+S^\top\bar{\Theta}+\bar{\Theta}^\top R\bar{\Theta}\big]X,X\big\rangle+2\big\langle \big[B^\top P+D^\top P(C+D\bar{\Theta})\big]X,v  \big\rangle \\
&\qquad +\int_\mathbf{E} \big\langle \big[ \textcolor{red}{2}G(e)^\top P\big(F(e)+G(e)\bar{\Theta}\big)X+G(e)^\top PG(e)v\big],v\big\rangle \pi(de)\\
&\qquad\qquad\quad +\big\langle D^\top PDv,v\big\rangle\bigg\}ds+[\cdots]dW+\int_\mathbf{E}[\cdots]\tilde{\mathnormal{N}}(deds).
\end{aligned}
\end{equation*}
Thus, we have
\begin{equation*}
\begin{aligned}
&\mathbb{E}\big\langle HX(T),X(T) \big\rangle-\mathbb{E}\big\langle P(t)x,x \big\rangle\\
&=\mathbb{E}\int_t^T\biggl\{ -\Big\langle \big[Q+\bar{\Theta}^\top S+S^\top\bar{\Theta}+\bar{\Theta}^\top R\bar{\Theta}\big]X,X \Big\rangle
+\Big\langle \big[D^\top PD+\int_\mathbf{E}G(e)^\top PG(e)\pi(de)\big] v,v \Big\rangle  \\
&\qquad\qquad +2\Big\langle \big[B^\top P+D^\top P(C+D\bar{\Theta})+\int_\mathbf{E}G(e)^\top P\big(F(e)+G(e)\bar{\Theta}\big)\pi(de)\big]X,v \Big\rangle\biggr\}ds.
\end{aligned}
\end{equation*}
Putting the above equation into the cost functional, we have
\begin{equation*}
\begin{aligned}
\tilde{J}(t,x;v(\cdot))&=\mathbb{E}\big\langle P(t)x,x \big\rangle+\mathbb{E}\int_t^T\biggl\{
\Big\langle \big[ R+D^\top PD+\int_\mathbf{E} G(e)^\top PG(e)\pi(de)\big]v,v \Big\rangle\\
&\qquad+2\Big\langle \big[ B^\top P+D^\top PC+S+\int_\mathbf{E}G(e)^\top PF(e)\pi(de)\\
&\qquad+\big(R+D^\top PD+\int_\mathbf{E}G(e)^\top PG(e)\pi(de)\big)\bar{\Theta} \big]X,v \Big\rangle\biggr\}ds.
\end{aligned}
\end{equation*}
From the previous analysis, we can know that $\bar{v}(\cdot)=0$ is the open-loop optimal control for Problem (SLQP)$^0_{\bar{\Theta}}$, i.e.,
\begin{equation*}
\tilde{J}(t,x;v(\cdot)) \geqslant \tilde{J}(t,x;0)=\mathbb{E}\big\langle P(t)x,x \big\rangle,
\end{equation*}
then
\begin{equation*}
\begin{aligned}
0&\leqslant\tilde{J}(t,x;v(\cdot))-\mathbb{E}\big\langle P(t)x,x \big\rangle
=\mathbb{E}\int_{t}^{T}\biggl\{\Big\langle \big[ R+D^\top PD+\int_\mathbf{E} G(e)^\top PG(e)\pi(de)\big]v,v \Big\rangle\\
&\qquad +2\Big\langle \big[ B^\top P+D^\top PC+S+\int_\mathbf{E}G(e)^\top PF(e)\pi(de)\\
&\qquad +\big(R+D^\top PD+\int_\mathbf{E}G(e)^\top PG(e)\pi(de)\big)\bar{\Theta} \big]X,v \Big\rangle\biggr\}ds.
\end{aligned}
\end{equation*}

For simplicity, we give the following notations:
\begin{equation*}\left\{
\begin{aligned}
\mathcal{R}&:=R+D^\top PD+\int_\mathbf{E} G(e)^\top PG(e)\pi(de),\\
\mathcal{L}&:=B^\top P+D^\top PC+S+\int_\mathbf{E}G(e)^\top PF(e)\pi(de),
\end{aligned}\right.
\end{equation*}
thus,
\begin{equation*}
\mathbb{E}\int_t^T\Big[\big\langle \mathcal{R}v,v \big\rangle+2\big\langle \big(\mathcal{L}+\mathcal{R}\bar{\Theta}\big)X,v \big\rangle\Big]ds \geqslant 0, \,\,\forall v(\cdot) \in \mathcal{U}[t,T].
\end{equation*}
Choose the initial pair $(t,x)= (0,0)$ and $v(\cdot)=v_0I_{[r,r+h]}(\cdot)$, $0\leqslant r\textless r+h \leqslant T$, $v_0 \in \mathbb{R}^m$. In this time, $v(\cdot)$ is a deterministic function, thus
\begin{equation*}
\begin{aligned}
&\mathbb{E}\int_0^T\Big[\big\langle \mathcal{R}v,v \big\rangle+2\big\langle \big(\mathcal{L}+\mathcal{R}\bar{\Theta}\big)X,v \big\rangle\Big]ds
=\int_0^T\Big[\big\langle \mathcal{R}v,v \big\rangle+2\big\langle \big(\mathcal{L}+\mathcal{R}\bar{\Theta}\big)\mathbb{E}X,v \big\rangle\Big]ds\\
&=\int_r^{r+h}\Big[\big\langle \mathcal{R}v_0,v_0 \big\rangle+2\big\langle \big(\mathcal{L}+\mathcal{R}\bar{\Theta}\big)\mathbb{E}X,v_0 \big\rangle\Big]ds,
\end{aligned}
\end{equation*}
where $\mathbb{E}X(\cdot)$ is solution to the following {\it ordinary differential equation} (ODE):
\begin{equation*}\left\{
\begin{aligned}
d\mathbb{E}X(s)&=\Big\{\big[A(s)+B(s)\bar{\Theta}(s)\big]\mathbb{E}X(s)+B(s)v_0I_{[r,r+h]}(s)\Big\}ds,\quad s \in [0,T],\\
\mathbb{E}X(0)&=0.
\end{aligned}\right.
\end{equation*}
Then
\begin{equation*}
\mathbb{E}X(t)=\int_{0}^{t}B(s)v_0I_{[r,r+h]}(s)e^{-\int_{0}^{s}\big(A(r)+B(r)\bar{\Theta}(r)\big)dr}ds\cdot e^{\int_{0}^{t}\big(A(r)+B(r)\bar{\Theta}(r)\big)dr}, \,\,t \in [0,T].
\end{equation*}
It is easy to see that
\begin{equation*}
\mathbb{E}X(t) \to 0,\,\,\mbox{as}\,\,h \to 0,\,\,t \in [0,T],
\end{equation*}
thus
\begin{equation*}
\lim_{h \to 0}\frac{1}{h}\int_r^{r+h}\Big[\big\langle \mathcal{R}(s)v_0,v_0 \big\rangle
+2\big\langle \big(\mathcal{L}(s)+\mathcal{R}(s)\bar{\Theta}(s)\big)\mathbb{E}X(s),v_0 \big\rangle\Big]ds \geqslant 0.
\end{equation*}
We can know
\begin{equation*}
\big\langle \mathcal{R}(r)v_0,v_0 \big\rangle \geqslant 0,\,\,r \in [0,T],\,\,\forall v_0 \in \mathbb{R}^m \quad \Rightarrow \quad \mathcal{R}(r) \geqslant 0,\,\,r \in [0,T].
\end{equation*}
Thus, the first inequality in (\ref{regularity}) is obtained and now let us prove the second equation. We take the initial pair $(t,x) \in [0,T) \times \mathbb{R}^n$ and $v(\cdot)=\frac{v_0}{n}I_{[r,r+h]}(\cdot)$, $v_0 \in \mathbb{R}^m$, $t\leqslant r\textless r+h \leqslant T$. In this time,
\begin{equation*}
\begin{aligned}
&\mathbb{E}\int_t^T\Big[\big\langle \mathcal{R}v,v \big\rangle+2\big\langle \big(\mathcal{L}+\mathcal{R}\bar{\Theta}\big)X,v \big\rangle\Big]ds
=\mathbb{E}\int_r^{r+h}\biggl[\big\langle \mathcal{R}\frac{v_0}{n},\frac{v_0}{n} \big\rangle+2\big\langle \big(\mathcal{L}+\mathcal{R}\bar{\Theta}\big)X,\frac{v_0}{n} \big\rangle\biggr]ds\\
&=\int_r^{r+h}\biggl[\frac{1}{n^2}\big\langle \mathcal{R}v_0,v_0 \big\rangle+\frac{2}{n}\big\langle \big(\mathcal{L}+\mathcal{R}\bar{\Theta}\big)\mathbb{E}X,v_0 \big\rangle\biggr]ds,
\end{aligned}
\end{equation*}
where $\mathbb{E}X(\cdot)$ is solution to the following ODE:
\begin{equation*}\left\{
\begin{aligned}
d\mathbb{E}X(s)&=\bigg\{\big[A(s)+B(s)\bar{\Theta}(s)\big]\mathbb{E}X(s)+\frac{B(s)v_0}{n}I_{[r,r+h]}(s)\bigg\}ds,\quad s \in [t,T],\\
\mathbb{E}X(t)&=\mathbb{E}x.
\end{aligned}\right.
\end{equation*}
As we know, this is an inhomogeneous linear ODE, and we get for $s \in [t,T]$,
\begin{equation*}
\mathbb{E}X(s)=\bigg[\mathbb{E}x+\int_{t}^{s}B(y)\frac{v_0}{n}I_{[r,r+h]}(y)e^{-\int_{t}^{y}\big(A(z)+B(z)\bar{\Theta}(z)\big)dz}dy\bigg]\cdot e^{\int_{t}^{s}\big(A(z)+B(z)\bar{\Theta}(z)\big)dz}.
\end{equation*}
It is easy to see
\begin{equation*}
\mathbb{E}X(s) \to \mathbb{E}x\cdot e^{\int_{t}^{s}\big(A(z)+B(z)\bar{\Theta}(z)\big)dz} \triangleq \mathbb{E}x\cdot\Psi(s) ,\,\,\mbox{as}\,\,n \to \infty.
\end{equation*}
Thus
\begin{equation*}
\begin{aligned}
&\lim_{n \to \infty}\int_{r}^{r+h}n\cdot\biggl[\frac{1}{n^2}\big\langle \mathcal{R}v_0,v_0 \big\rangle+\frac{2}{n}\big\langle \big(\mathcal{L}(s)+\mathcal{R}(s)\bar{\Theta}(s)\big)\mathbb{E}X(s),v_0 \big\rangle\biggr]ds\\
&= \lim_{n \to \infty}\int_{r}^{r+h}\cdot\biggl[\frac{1}{n}\big\langle \mathcal{R}v_0,v_0 \big\rangle+2\big\langle \big(\mathcal{L}(s)+\mathcal{R}(s)\bar{\Theta}(s)\big)\mathbb{E}X(s),v_0 \big\rangle\biggr]ds\\
&= 2\int_{r}^{r+h}\big\langle \big(\mathcal{L}(s)+\mathcal{R}(s)\bar{\Theta}(s)\big)\Psi(s)\mathbb{E}x,v_0 \big\rangle ds \geqslant 0.
\end{aligned}
\end{equation*}
Since $v_0 \in \mathbb{R}^m$ is arbitrary,
\begin{equation*}
\int_{r}^{r+h}\big\langle \big(\mathcal{L}(s)+\mathcal{R}(s)\bar{\Theta}(s)\big)\Psi(s)\mathbb{E}x,v_0 \big\rangle ds=0.
\end{equation*}
Divide both sides by $h$ and let $h \to 0$, we have
\begin{equation*}
\big\langle \big[\mathcal{L}(r)+\mathcal{R}(r)\bar{\Theta}(r)\big]\Psi(r)\mathbb{E}x,v_0 \big\rangle =0,\,\,\forall r \in [t,T].
\end{equation*}
Since $v_0 \in \mathbb{R}^m$ is arbitrary,
\begin{equation*}
\mathcal{L}(r)+\mathcal{R}(r)\bar{\Theta}(r)=0,\,\,\forall r \in [t,T].
\end{equation*}
This is the second equation in (\ref{regularity}). The proof is complete. $\qquad\Box$

The following theorem is the main result in this paper, which characterizes the closed-loop solvability of Problem (SLQP).

\begin{mythm}\label{clos}
Let (H1)-(H2) hold. Problem (SLQP) admits an optimal closed-loop strategy $(\bar{\Theta}(\cdot),\bar{v}(\cdot)) \in \mathcal{Q}[t,T] \times \mathcal{U}[t,T]$ if and only if the following RIDE:
\begin{equation}\label{Riccati}\left\{
\begin{aligned}
&0=\dot{P}+PA+A^\top P+C^\top PC+\int_\mathbf{E}F(e)^\top PF(e)\pi(de)+Q\\
&\qquad -\bigg(B^\top P +D^\top PC+\int_\mathbf{E}G(e)^\top PF(e)\pi(de)+S\bigg)^\top\bigg(R+D^\top PD\\
&\qquad\quad  +\int_\mathbf{E}G(e)^\top PG(e)\pi(de)\bigg)^\dagger\bigg(B^\top P +D^\top PC+\int_\mathbf{E}G(e)^\top PF(e)\pi(de)+S\bigg),\, \textcolor{red}{s \in [t,T]},\\
&\mathcal{R}\bigg(B^\top P +D^\top PC+\int_\mathbf{E}G(e)^\top PF(e)\pi(de)+S\bigg)\\
&\quad \subseteq \mathcal{R}\bigg(R+D^\top PD+\int_\mathbf{E}G(e)^\top PG(e)\pi(de)\bigg),\\
&R+D^\top PD+\int_\mathbf{E}G(e)^\top PG(e)\pi(de) \geqslant 0,\quad P(T)=H,\\
\end{aligned}\right.
\end{equation}	
admits a solution $P(\cdot) \in C([t,T];\mathbb{S}^n)$ such that
\begin{equation}\label{range1}
\begin{aligned}
&\bigg(R+D^\top PD+\int_\mathbf{E}G(e)^\top PG(e)\pi(de)\bigg)^\dagger\bigg(B^\top P +D^\top PC+\int_\mathbf{E}G(e)^\top PF(e)\pi(de)+S\bigg)\\
&\quad \in L^2(t,T;\mathbb{R}^{m \times n}),
\end{aligned}\end{equation}
and the following BSDEP:
{\small\begin{equation}\label{follower BSDE}\left\{
\begin{aligned}
d\eta&=-\Bigg\{ \bigg[A-B\bigg(R+D^\top PD+\int_\mathbf{E}G(e)^\top PG(e)\pi(de)\bigg)^\dagger\\
&\qquad\qquad \times\bigg(B^\top P +D^\top PC+\int_\mathbf{E}G(e)^\top PF(e)\pi(de)+S\bigg)\bigg]^\top\eta\\
&\qquad+\bigg[C-D\bigg(R+D^\top PD+\int_\mathbf{E}G(e)^\top PG(e)\pi(de)\bigg)^{\dagger}\\
&\qquad\qquad \times\bigg(B^\top P +D^\top PC+\int_\mathbf{E}G(e)^\top PF(e)\pi(de)+S\bigg)\bigg]^\top(\zeta+P\sigma)\\
&\qquad+\int_\mathbf{E}\bigg[F(e)-G(e)\bigg(R+D^\top PD+\int_\mathbf{E}G(e)^\top PG(e)\pi(de)\bigg)^{\dagger}\\
&\qquad\qquad \times\bigg(B^\top P +D^\top PC+\int_\mathbf{E}G(e)^\top PF(e)\pi(de)+S\bigg)\bigg]^\top\big(\psi(e)+Pf(e)\big)\pi(de)\\
&\qquad +Pb+q-\bigg(B^\top P +D^\top PC+\int_\mathbf{E}G(e)^\top PF(e)\pi(de)+S\bigg)^\top\\
&\qquad\,\,\,\,\, \times\bigg(R+D^\top PD+\int_\mathbf{E}G(e)^\top PG(e)\pi(de)\big)^{\dagger}\rho\Bigg\}ds+\zeta dW+\int_\mathbf{E}\psi(e)\tilde{\mathnormal{N}}(deds),\,s \in [t,T],\\
&\hspace{-5mm}B^\top\eta+D^\top(\zeta+P\sigma)+\int_\mathbf{E}G(e)^\top\big(\psi(e)+Pf(e)\big)\pi(de)+\rho\\
& \in \mathcal{R}\big(R+D^\top PD+\int_\mathbf{E}G(e)^\top PG(e)\pi(de)\big),\,\,\, a.e.,\,\,\,\mathbb{P}\mbox{-}a.s.,\quad \eta(T)=g,
\end{aligned}\right.
\end{equation}}	
admits a unique solution $(\eta(\cdot),\zeta(\cdot),\psi(\cdot))$ which satisfies
\begin{equation}\label{range2}
\begin{aligned}
&\bigg(R+D^\top PD+\int_\mathbf{E}G(e)^\top PG(e)\pi(de)\bigg)^\dagger\\
&\quad \times\bigg(B^\top\eta+D^\top(\zeta+P\sigma)+\int_\mathbf{E}G(e)^\top\big(\psi(e)+Pf(e)\big)\pi(de)+\rho\bigg) \in L^2_{\mathbb{F},p}(t,T;\mathbb{R}^m).
\end{aligned}
\end{equation}
In this case, the closed-loop optimal strategy $(\bar{\Theta}(\cdot),\bar{v}(\cdot))$ of Problem (SLQP) admits the following representation:
\begin{equation}\label{follower closed-loop optimal}\left\{
\begin{aligned}
\bar{\Theta}&=-\bigg(R+D^\top PD+\int_\mathbf{E}G(e)^\top PG(e)\pi(de)\bigg)^\dagger\\
&\qquad \times\bigg(B^\top P +D^\top PC+\int_\mathbf{E}G(e)^\top PF(e)\pi(de)+S\bigg)\\
&\quad+\bigg[I-\bigg(R+D^\top PD+\int_\mathbf{E}G(e)^\top PG(e)\pi(de)\bigg)^\dagger\\
&\qquad \times\bigg(R+D^\top PD+\int_\mathbf{E}G(e)^\top PG(e)\pi(de)\bigg)\bigg]\theta,\\
\bar{v}&=-\bigg(R+D^\top PD+\int_\mathbf{E}G(e)^\top PG(e)\pi(de)\bigg)^\dagger\\
&\qquad \times\bigg[B^\top\eta+D^\top(\zeta+P\sigma)+\int_\mathbf{E}G(e)^\top\big(\psi(e)+Pf(e)\big)\pi(de)+\rho\bigg]\\
&\quad+\big[I-\bigg(R+D^\top PD+\int_\mathbf{E}G(e)^\top PG(e)\pi(de)\bigg)^\dagger\\
&\qquad \times\bigg(R+D^\top PD+\int_\mathbf{E}G(e)^\top PG(e)\pi(de)\bigg)\bigg]v,
\end{aligned}\right.
\end{equation}
for some $\theta(\cdot) \in \mathcal{Q}[t,T]$, $v(\cdot)\in \mathcal{U}[t,T]$. Further, the value function $V(\cdot,\cdot)$ is given by
\begin{equation}\label{value}
\begin{split}
V(t,x)&\equiv\inf_{u(\cdot)\in\,\mathcal{U}[t,T]}J(t,x;u(\cdot))\\
&=\mathbb{E}\biggl\{ \big\langle P(t)x,x\big\rangle +2\big\langle \eta(t),x\big\rangle+\int_t^T\bigg[ 2\big\langle \eta,b\big\rangle+2\big\langle \zeta,\sigma\big\rangle\\
&\qquad\quad +2\int_\mathbf{E}\big\langle \psi(e),f(e)\big\rangle\pi(de)+\big\langle P\sigma,\sigma \big\rangle+\int_\mathbf{E}\big\langle Pf(e),f(e)\big\rangle\pi(de)\\
&\qquad\quad -\bigg|\bigg[\bigg(R+D^\top PD+\int_\mathbf{E}G(e)^\top PG(e)\pi(de)\bigg)^\dagger\bigg]^{\frac{1}{2}}\\
&\qquad\qquad \times\bigg[B^\top\eta+D^\top(\zeta+P\sigma)+\int_\mathbf{E}G(e)^\top\big(\psi(e)+Pf(e)\big)\pi(de)+\rho\bigg]\bigg|^2\bigg]ds \biggr\}.
\end{split}
\end{equation}
\end{mythm}

\textit{Proof.} We first prove the necessity. Let $(\bar{\Theta}(\cdot),\bar{v}(\cdot)) \in \mathcal{Q}[t,T] \times \mathcal{U}[t,T]$ be a closed-loop optimal strategy of Problem (SLQP) over $[t,T]$. The second equation of (\ref{regularity}) implies
\begin{equation*}
\mathcal{R}\bigg(B^\top P +D^\top PC+\int_\mathbf{E}G(e)^\top PF(e)\pi(de)+S\bigg) \subseteq \mathcal{R}\bigg(R+D^\top PD+\int_\mathbf{E}G(e)^\top PG(e)\pi(de)\bigg),\,\,a.e..
\end{equation*}
Denoting $\mathcal{R}:= R+D^\top PD+\int_\mathbf{E} G(e)^\top PG(e)\pi(de)$, since
\begin{equation*}
\mathcal{R}^\dagger\bigg(B^\top P +D^\top PC+\int_\mathbf{E}G(e)^\top PF(e)\pi(de)+S\bigg)=-\mathcal{R}^\dagger\mathcal{R}\bar{\Theta},
\end{equation*}
and $\mathcal{R}^\dagger\mathcal{R}$ is an orthogonal projection, we see that (\ref{range1}) holds and
\begin{equation*}
\bar{\Theta}=-\mathcal{R}^\dagger\bigg(B^\top P +D^\top PC+\int_\mathbf{E}G(e)^\top PF(e)\pi(de)+S\bigg)+\big(I-\mathcal{R}^\dagger\mathcal{R}\big)\theta,
\end{equation*}
for some $\theta(\cdot) \in \mathcal{Q}[t,T]$. Consequently,
\begin{equation*}
\begin{aligned}
&\bigg(PB +C^\top PD+\int_\mathbf{E}F(e)^\top PG(e)\pi(de)+S^\top\bigg)\bar{\Theta}\\
=&-\bar{\Theta}^\top\mathcal{R}\bar{\Theta}=\bar{\Theta}^\top\mathcal{R}\mathcal{R}^\dagger\bigg(B^\top P +D^\top PC+\int_\mathbf{E}G(e)^\top PF(e)\pi(de)+S\bigg)\\
=&-\bigg(B^\top P +D^\top PC+\int_\mathbf{E}G(e)^\top PF(e)\pi(de)+S\bigg)^\top\mathcal{R}^\dagger\\
&\quad \times\bigg(B^\top P +D^\top PC+\int_\mathbf{E}G(e)^\top PF(e)\pi(de)+S\bigg).
\end{aligned}
\end{equation*}
Plug the above into the Lyapunov equation (\ref{Lyapunov}), we can obtain the RIDE in (\ref{Riccati}).

To determine $\bar{v}(\cdot)$, we define
\begin{equation*}\left\{
\begin{aligned}
\eta&=\bar{Y}-P\bar{X},\\
\zeta&=\bar{Z}-P(C+D\bar{\Theta})\bar{X}-PD\bar{v}-P\sigma,\\
\psi(e)&=\bar{K}(e)-P\big(F(e)+G(e)\bar{\Theta}\big)\bar{X}-PG(e)\bar{v}-Pf(e),
\end{aligned}\right.
\end{equation*}
where $(\bar{X}(\cdot),\bar{Y}(\cdot),\bar{Z}(\cdot),\bar{K}(\cdot,\cdot))$ is the adapted solution to the FBSDEP (\ref{Th3.1-os}). Then
\begin{equation*}
\begin{aligned}
d\eta&=d\bar{Y}-\dot{P}\bar{X}ds-Pd\bar{X}\\
&=\bigg\{-A^\top\bar{Y}-C^\top\bar{Z}-\int_\mathbf{E}F(e)^\top\bar{K}(e)\pi(de)-(Q+S^\top\bar{\Theta})\bar{X}-S^\top\bar{v}-q\\
&\qquad +A^\top P\bar{X}+PA\bar{X}+C^\top PC\bar{X}+\bigg(PB+C^\top PD+\int_\mathbf{E}F(e)^\top PG(e)\pi(de)\\
&\qquad +S^\top\bigg)\bar{\Theta}\bar{X}+Q\bar{X}+\int_\mathbf{E}F(e)^\top PF(e)\bar{X}\pi(de)-P(A+B\bar{\Theta})\bar{X}-PB\bar{v}-Pb\bigg\}ds\\
&\quad +\big[\bar{Z}-P(C+D\bar{\Theta})\bar{X}-PD\bar{v}-P\sigma\big]dW\\
&\quad +\int_\mathbf{E}\big[\bar{K}(e)-P\big(F(e)+G(e)\bar{\Theta}\big)\bar{X}_{-}-PG(e)\bar{v}-Pf(e)\big]\tilde{\mathnormal{N}}(deds)\\
&=-\bigg\{A^\top(P\bar{X}+\eta)+C^\top\big[\zeta+P(C+D\bar{\Theta})\bar{X}+PD\bar{v}+P\sigma\big]\\
&\qquad\ +\int_\mathbf{E}F(e)^\top\big[\psi(e)+P\big(F(e)+G(e)\bar{\Theta}\big)\bar{X}_{-}+PG(e)\bar{v}+Pf(e)\big]\pi(de)\\
\end{aligned}
\end{equation*}
\begin{equation*}
\begin{aligned}
&\qquad\ +(Q+S^\top\bar{\Theta})\bar{X}+S^\top\bar{v}+q-A^\top P\bar{X}-PA\bar{X}-C^\top PC\bar{X}\\
&\qquad\ -\bigg(PB+C^\top PD+\int_\mathbf{E}F(e)^\top PG(e)\pi(de)+S^\top\bigg)\bar{\Theta}\bar{X}-Q\bar{X}\\
&\qquad\ -\int_\mathbf{E}F(e)^\top PF(e)\bar{X}\pi(de)+P(A+B\bar{\Theta})\bar{X}+PB\bar{v}+Pb\bigg\}ds\\
&\qquad\ +\zeta dW+\int_\mathbf{E}\psi(e)\tilde{\mathnormal{N}}(deds)\\
&=-\bigg\{A^\top\eta+C^\top\zeta+\int_\mathbf{E}F(e)^\top\psi(e)\pi(de)+\bigg(PB+C^\top PD+\int_\mathbf{E}F(e)^\top PG(e)\pi(de)\\
&\qquad\ +S^\top\bigg)\bar{v}+Pb+C^\top P\sigma+\int_\mathbf{E}F(e)^\top Pf(e)\pi(de)+q\bigg\}ds+\zeta dW+\int_\mathbf{E}\psi(e)\tilde{\mathnormal{N}}(deds).
\end{aligned}
\end{equation*}

According to the stationarity condition (\ref{Th3.1-sc}), we have
\begin{equation}\label{v_1 condition}
\begin{aligned}
0&=B^\top \bar{Y}+D^\top\bar{Z}+\int_\mathbf{E}G(e)^\top\bar{K}(e)\pi(de)+(S+R\bar{\Theta})\bar{X}+R\bar{v}+\rho\\
&=B^\top (P\bar{X}+\eta)+D^\top\big[\zeta+P(C+D\bar{\Theta})\bar{X}+PD\bar{v}+P\sigma\big]+(S+R\bar{\Theta})\bar{X}+R\bar{v}+\rho\\
&\qquad +\int_\mathbf{E}G(e)^\top \big[\psi(e)+P\big(F(e)+G(e)\bar{\Theta}\big)\bar{X}+PG(e)\bar{v}+Pf(e)\big] \pi(de)\\
&=B^\top\eta+D^\top\zeta+\int_\mathbf{E}G(e)^\top \psi(e)\pi(de)+D^\top P\sigma+\rho+\int_\mathbf{E}G(e)^\top Pf(e)\pi(de)\\
&\quad +\bigg(R+D^\top PD+\int_\mathbf{E}G(e)^\top PG(e)\pi(de)\bigg)\bar{v}+\bigg[B^\top P+D^\top PC\\
&\quad +\int_\mathbf{E}G(e)^\top PF(e)\pi(de)+S+\bigg(R+D^\top PD+\int_\mathbf{E}G(e)^\top PG(e)\pi(de)\bigg)\bar{\Theta}\bigg]\bar{X}\\
&=B^\top\eta+D^\top\zeta+\int_\mathbf{E}G(e)^\top \psi(e)\pi(de)+D^\top P\sigma+\rho+\int_\mathbf{E}G(e)^\top Pf(e)\pi(de)\\
&\quad +\bigg(R+D^\top PD+\int_\mathbf{E}G(e)^\top PG(e)\pi(de)\bigg)\bar{v}.
\end{aligned}
\end{equation}
Hence,
\begin{equation}
\begin{aligned}
&B^\top\eta+D^\top\zeta+\int_\mathbf{E}G(e)^\top \psi(e)\pi(de)+D^\top P\sigma+\rho+\int_\mathbf{E}G(e)^\top Pf(e)\pi(de) \\
&\qquad \in \mathcal{R}\bigg(R+D^\top PD+\int_\mathbf{E}G(e)^\top PG(e)\pi(de)\bigg),\quad a.e.,\ \mathbb{P}\mbox{-}a.s..
\end{aligned}
\end{equation}
Since
\begin{equation}\nonumber
\mathcal{R}^\dagger\bigg[B^\top\eta+D^\top\zeta+\int_\mathbf{E}G(e)^\top \psi(e)\pi(de)+D^\top P\sigma+\rho+\int_\mathbf{E}G(e)^\top Pf(e)\pi(de)\bigg]=-\mathcal{R}^\dagger\mathcal{R}\bar{v}
\end{equation}
and $\mathcal{R}^\dagger\mathcal{R}$ is an orthogonal projection, we see that (\ref{range2}) holds and
\begin{equation}\nonumber
\bar{v}=-\mathcal{R}^\dagger\bigg[B^\top\eta+D^\top(\zeta+P\sigma)+\int_\mathbf{E}G(e)^\top \big(\psi(e)+Pf(e)\big)\pi(de)+\rho\bigg]+(I-\mathcal{R}^\dagger\mathcal{R})v,
\end{equation}
for some $v(\cdot) \in L^2_{\mathbb{F},p}(t,T;\mathbb{R}^m)$.
Consequently,
\begin{equation*}
\begin{aligned}
&\bigg[PB+C^\top PD+\int_\mathbf{E}F(e)^\top PG(e)\pi(de)+S^\top\bigg]\bar{v}\\
&=-\bigg[PB+C^\top PD+\int_\mathbf{E}F(e)^\top PG(e)\pi(de)+S^\top\bigg]\\
&\qquad \times \mathcal{R}^\dagger\bigg[B^\top\eta+D^\top(\zeta+P\sigma)+\int_\mathbf{E}G(e)^\top \big(\psi(e)+Pf(e)\big)\pi(de)+\rho\bigg]\\
&\quad +\bigg[PB+C^\top PD+\int_\mathbf{E}F(e)^\top PG(e)\pi(de)+S^\top\bigg](I-\mathcal{R}^\dagger\mathcal{R})v\\
&=-\bigg[PB+C^\top PD+\int_\mathbf{E}F(e)^\top PG(e)\pi(de)+S^\top\bigg]\\
&\qquad \times \mathcal{R}^\dagger\bigg[B^\top\eta+D^\top(\zeta+P\sigma)+\int_\mathbf{E}G(e)^\top \big(\psi(e)+Pf(e)\big)\pi(de)+\rho\bigg]\\
&\quad -\bar{\Theta}^\top\mathcal{R}(I-\mathcal{R}^\dagger\mathcal{R})v\\
&=-\big[PB+C^\top PD+\int_\mathbf{E}F(e)^\top PG(e)\pi(de)+S^\top\big]\\
&\qquad \times\mathcal{R}^\dagger\bigg[B^\top\eta+D^\top(\zeta+P\sigma)+\int_\mathbf{E}G(e)^\top \big(\psi(e)+Pf(e)\big)\pi(de)+\rho\bigg].
\end{aligned}
\end{equation*}
Therefore, $(\eta(\cdot),\zeta(\cdot),\psi(\cdot,\cdot))$ is the unique solution to the following BSDEP:
\begin{equation*}\nonumber\left\{
\begin{aligned}
d\eta&=-\bigg\{A^\top\eta+C^\top\zeta+\int_\mathbf{E}F(e)^\top\psi(e)\pi(de)\\
&\qquad +\bigg[PB+C^\top PD+\int_\mathbf{E}F(e)^\top PG(e)\pi(de)+S^\top\bigg]\bar{v}\\
&\qquad +Pb+C^\top P\sigma+\int_\mathbf{E}F(e)^\top Pf(e)\pi(de)+q\bigg\}ds+\zeta dW+\int_\mathbf{E}\psi(e)\tilde{\mathnormal{N}}(deds)\\
&=-\bigg\{A^\top\eta+C^\top\zeta+\int_\mathbf{E}F(e)^\top\psi(e)\pi(de)+Pb+C^\top P\sigma+\int_\mathbf{E}F(e)^\top Pf(e)\pi(de)+q\\
&\qquad\quad -\bigg[PB+C^\top PD+\int_\mathbf{E}F(e)^\top PG(e)\pi(de)+S^\top\bigg]\\
&\qquad\qquad \times \mathcal{R}^\dagger\bigg[B^\top\eta+D^\top(\zeta+P\sigma)+\int_\mathbf{E}G(e)^\top \big(\psi(e)+Pf(e)\big)\pi(de)+\rho\bigg]\bigg\}ds\\
&\qquad +\zeta dW+\int_\mathbf{E}\psi(e)\tilde{\mathnormal{N}}(deds),\quad s \in[t,T],\\
\eta(T)&=g.
\end{aligned}\right.
\end{equation*}

To prove the sufficiency, we take any $u(\cdot) \in \mathcal{U}[t,T]$, and let $X(\cdot)\equiv X(\cdot;t,x,u(\cdot))$ be the corresponding state process. Then, applying It\^o-Wentzell's formula to $\langle PX,X \rangle$ and $\langle \eta,X \rangle $, we have
\begin{equation}
\begin{aligned}
d\big\langle PX,X\rangle&=\bigg\langle -P(A+B\bar{\Theta})X-(A+B\bar{\Theta})^\top PX-(C+D\bar{\Theta})^\top P(C+D\bar{\Theta})X\\
&\qquad -\int_\mathbf{E}\big[(F(e)+G(e)\bar{\Theta})^\top P(F(e)+G(e)\bar{\Theta})X\big]\pi(de)\\
&\qquad -(Q+\bar{\Theta}^\top S+S^\top\bar{\Theta}+\bar{\Theta}^\top R\bar{\Theta})X,X\bigg\rangle ds\\
&\quad +2\big\langle PAX+PBu+Pb,X\big\rangle ds+ \big\langle P(CX+Du+\sigma),CX+Du+\sigma\big\rangle ds\\
&\quad +\int_\mathbf{E}\Big\langle P\big[F(e)X+G(e)u+f(e)\big],F(e)X+G(e)u+f(e)\Big\rangle \pi(de)ds\\
&\quad +[\cdots]dW+\int_\mathbf{E}[\cdots]\tilde{\mathnormal{N}}(deds)\\
&=\bigg\{-\bigg\langle \bigg[(PB+S^\top)\bar{\Theta}+\bar{\Theta}^\top (B^\top P+S)+(C+D\bar{\Theta})^\top P(C+D\bar{\Theta})\\
&\qquad +\int_\mathbf{E}\big[(F(e)+G(e)\bar{\Theta})^\top P(F(e)+G(e)\bar{\Theta})\big]\pi(de)+Q+\bar{\Theta}^\top R\bar{\Theta}\bigg]X,X\bigg\rangle\\
&\qquad +2\big\langle P(Bu+b),X\big\rangle+\big\langle P(CX+Du+\sigma),CX+Du+\sigma\big\rangle\\
&\qquad +\int_\mathbf{E}\Big\langle P\big[F(e)X+G(e)u+f(e)\big],F(e)X+G(e)u+f(e)\Big\rangle \pi(de)\bigg\}ds\\
&\quad +[\cdots]dW+\int_\mathbf{E}[\cdots]\tilde{\mathnormal{N}}(deds),\\
\end{aligned}
\end{equation}
and
\begin{equation*}
\begin{aligned}
d\big\langle \eta,X \big\rangle &=\biggl\{\bigg\langle -A^\top\eta-C^\top\zeta-\int_\mathbf{E}F(e)^\top\psi(e)\pi(de)\\
&\qquad\quad -\bigg[PB+C^\top PD+\int_\mathbf{E}F(e)^\top PG(e)\pi(de)+S^\top\bigg]\bar{v}\\
&\qquad\quad -Pb-C^\top P\sigma-\int_\mathbf{E}F(e)^\top Pf(e)\pi(de)-q,X \bigg\rangle+\big\langle \eta,AX+Bu+b\big\rangle\\
&\qquad +\big\langle \zeta, CX+Du+\sigma\big\rangle+\int_\mathbf{E}\big\langle \psi(e),F(e)X+G(e)u+f(e)\big\rangle\pi(de) \biggr\}ds\\
&\qquad+[\cdots]dW+\int_\mathbf{E}[\cdots]\tilde{\mathnormal{N}}(deds)\\
&=\biggl\{-\bigg\langle \bigg[PB+C^\top PD+\int_\mathbf{E}F(e)^\top PG(e)\pi(de)+S^\top\bigg]\bar{v}+Pb+C^\top P\sigma\\
&\qquad\qquad +\int_\mathbf{E}F(e)^\top Pf(e)\pi(de)+q,X \bigg\rangle+\big\langle \eta,Bu+b\big\rangle+\big\langle \zeta, Du+\sigma\big\rangle\\
&\qquad +\int_\mathbf{E}\big\langle \psi(e),G(e)u+f(e)\big\rangle\pi(de) \biggr\}ds+[\cdots]dW+\int_\mathbf{E}[\cdots]\tilde{\mathnormal{N}}(deds).\\
\end{aligned}
\end{equation*}
Thus,
\begin{equation*}
\begin{aligned}
&J(t,x;u(\cdot))=\mathbb{E}
\bigg\{\big\langle HX(T),X(T)\big\rangle+2\big\langle g,X(T)\big\rangle \\
&\qquad\qquad\qquad\qquad +\int_t^T \Big[\big\langle QX,X\big\rangle+2\big\langle SX,u\big\rangle+\big\langle Ru,u\big\rangle+2\big\langle q,X\big\rangle+2\big\langle \rho,u\big\rangle\Big] ds\bigg\}\\
&=\mathbb{E}\Bigg\{\big\langle P(t)x,x\big\rangle+2\big\langle \eta(t),x\big\rangle +\int_t^T \bigg[\big\langle QX,X\big\rangle+2\big\langle SX,u\big\rangle+\big\langle Ru,u\big\rangle+2\big\langle q,X\big\rangle\\
&\qquad\quad +2\big\langle \rho,u\big\rangle-2\big\langle B^\top PX,\bar{\Theta}X \big\rangle-\big\langle C^\top PCX,x \big\rangle-2\big\langle D^\top PCX,\bar{\Theta}X\big\rangle\\
&\qquad\quad -\big\langle D^\top PD\bar{\Theta}X,\bar{\Theta}X \big\rangle-\big\langle QX,X \big\rangle-2\big\langle SX,\bar{\Theta}X\big\rangle-\big\langle R\bar{\Theta}X,\bar{\Theta}X\big\rangle\\
&\qquad\quad +2\big\langle B^\top PX,u\big\rangle+2\big\langle Pb,X \big\rangle+\big\langle C^\top PCX,X \big\rangle+\big\langle D^\top PCX,u \big\rangle\\
&\qquad\quad +2\big\langle PCX,\sigma\big\rangle+\big\langle D^\top PXC,u \big\rangle+\big\langle D^\top PDu,u \big\rangle+2\big\langle PDu,\sigma\big\rangle\\
&\qquad\quad +\big\langle P\sigma,\sigma \big\rangle-2\bigg\langle \bigg[PB+C^\top PD+\int_\mathbf{E}F(e)^\top PG(e)\pi(de)+S^\top\bigg]\bar{v},X \bigg\rangle\\
&\qquad\quad -2\bigg\langle C^\top P\sigma+Pb+\int_\mathbf{E}F(e)^\top Pf(e)\pi(de)+q,X \bigg\rangle+2\big\langle B^\top\eta,u \big\rangle+2\big\langle \eta,b \big\rangle\\
&\qquad\quad +2\big\langle D^\top\zeta,u \big\rangle+2\big\langle \zeta,\sigma \big\rangle-\int_\mathbf{E}\big\langle F(e)^\top PF(e)X,X \big\rangle\pi(de)\\
&\qquad\quad -2\int_\mathbf{E}\big\langle G(e)^\top PF(e)X,\bar{\Theta}X \big\rangle\pi(de)-\int_\mathbf{E}\big\langle G(e)^\top PG(e)\bar{\Theta}X,\bar{\Theta}X \big\rangle\pi(de)\\
&\qquad\quad +\int_\mathbf{E}\big\langle F(e)^\top PF(e)X,X \big\rangle\pi(de)+2\int_\mathbf{E}\big\langle G(e)^\top PF(e)X,u \big\rangle\pi(de)\\
&\qquad\quad +2\int_\mathbf{E}\big\langle F(e)^\top Pf(e),X \big\rangle\pi(de)+\int_\mathbf{E}\big\langle G(e)^\top PG(e)u,u \big\rangle\pi(de)\\
&\quad\qquad +2\int_\mathbf{E}\big\langle G(e)^\top Pf(e),u \big\rangle\pi(de)+\int_\mathbf{E}\big\langle Pf(e),f(e) \big\rangle\pi(de)\\
&\qquad\quad +2\int_\mathbf{E}\big\langle G(e)^\top\psi(e),u \big\rangle\pi(de)+2\int_\mathbf{E}\big\langle \psi(e),f(e) \big\rangle\pi(de)\bigg] ds\Bigg\}\\
&=\mathbb{E}\Bigg\{\big\langle P(t)x,x\big\rangle+2\big\langle \eta(t),x\big\rangle+\int_t^T \bigg[\big\langle \mathcal{R}u,u\big\rangle-\big\langle \mathcal{R}\bar{\Theta}X,\bar{\Theta}X\big\rangle\\
&\qquad\quad +2\bigg\langle B^\top P+D^\top PC+\int_\mathbf{E}G(e)^\top PF(e)\pi(de)+S,u-\bar{\Theta}X-\bar{v} \bigg\rangle\\
&\qquad\quad +2\bigg\langle B^\top\eta+D^\top\zeta+\int_\mathbf{E}G(e)^\top \psi(e)\pi(de)+D^\top P\sigma+\rho\\
&\qquad\quad +\int_\mathbf{E}G(e)^\top Pf(e)\pi(de),u \bigg\rangle +\big\langle P\sigma,\sigma\big\rangle+\int_\mathbf{E}\big\langle Pf(e),f(e)\big\rangle \pi(de)\\
&\qquad\quad +2\big\langle \eta,b \big\rangle+2\big\langle \zeta,\sigma \big\rangle+2\int_\mathbf{E}\big\langle \psi(e),f(e)\big\rangle\pi(de)\bigg] ds\Bigg\}\\
&=\mathbb{E}\Bigg\{\big\langle P(t)x,x\big\rangle+2\big\langle \eta(t),x\big\rangle+\int_t^T\big\langle \mathcal{R}\big(u-\bar{\Theta}X-\bar{v}\big),u-\bar{\Theta}X-\bar{v}\big\rangle ds\\
\end{aligned}
\end{equation*}
\begin{equation*}
\begin{aligned}
&\qquad+\int_t^T\bigg[\big\langle P\sigma,\sigma\big\rangle+\int_\mathbf{E}\big\langle Pf(e),f(e)\big\rangle \pi(de)+2\big\langle \eta,b \big\rangle+2\big\langle \zeta,\sigma \big\rangle
+2\int_\mathbf{E}\big\langle \psi(e),f(e)\big\rangle\pi(de)\\
&\qquad\qquad -\bigg\langle \mathcal{R}^\dagger\bigg[B^\top\eta+D^\top\zeta+\int_\mathbf{E}G(e)^\top \psi(e)\pi(de)+D^\top P\sigma+\rho+\int_\mathbf{E}G(e)^\top Pf(e)\pi(de)\bigg],\\
&\qquad\qquad\qquad\ B^\top\eta+D^\top\zeta+\int_\mathbf{E}G(e)^\top \psi(e)\pi(de)+D^\top P\sigma+\rho+\int_\mathbf{E}G(e)^\top Pf(e)\pi(de) \bigg\rangle\bigg] ds\Bigg\}\\
&=J(t,x;\bar{\Theta}(\cdot)\bar{X}(\cdot)+\bar{v}(\cdot))+\int_t^T\big\langle \mathcal{R}\big(u-\bar{\Theta}X-\bar{v}\big),u-\bar{\Theta}X-\bar{v}\big\rangle ds.
\end{aligned}
\end{equation*}
Hence,
\begin{equation}\nonumber
J(t,x;\bar{\Theta}(\cdot)\bar{X}(\cdot)+\bar{v}(\cdot)) \leqslant J(t,x;u(\cdot)),\quad \forall u(\cdot) \in \mathcal{U}[t,T],
\end{equation}
if and only if
\begin{equation}\nonumber
R+D^\top PD+\int_\mathbf{E} G(e)^\top PG(e)\pi(de) \geqslant 0,\quad a.e..
\end{equation}
This completes the proof. $\qquad\Box$

\section{Concluding remarks}

In this paper, we have investigated the closed-loop solvability of Problem (SLQP): the stochastic linear-quadratic optimal control problem with Poisson jumps.
We transform the initial problem into a new and simple problem and obtain a Riccati integral-differential equation firstly. Then we get the characterization of its closed-loop solvability, i.e., the solvabilities of the RIDE and a BSDEP (Theorem \ref{clos}).

Motivated by \cite{SY2018,SunXiongYong2021,SY2014,SY2019,LSY2021}, characterization of the closed-loop solvability for SLQP optimal control problem in infinite horizon, with random coefficients, and of the closed-loop saddle points, Nash equilibria for SLQ zero-sum, nonzero-sum differential games, respectively, are our future research topic. The closed-loop solvability for Stackelberg stochastic LQ differential game is recently studied by Li and Shi \cite{LiShi2021}, and we are also interested in that for Stackelberg stochastic LQ differential game of mean-field's type. We hope to report some results relevant to the above-mentioned problems in our future publications.


\begin{thebibliography}{0}

\bibitem{AMZ2001} M. Ait Rami, J. B. Moore and X. Y. Zhou, Indefinite stochastic linear quadratic control and generalized differential Riccati equation, \emph{SIAM J. Control Optim.}, {\bf 40}(4), 1296-1311, 2001.
	
\bibitem{Athans1968} H. Athans, The matrix minmum principle, \emph{Inform. Control}, {\bf 11}, 592-606, 1968.

\bibitem{AndersonMoore1989}B. D. O. Anderson, J. B. Moore, \emph{Optimal Control: Linear Quadratic Methods}, Prentice Hall, Englewood Cliffs, New Jersey, 1989.

\bibitem{Bensoussan1981} A. Bensoussan, Lectures on stochastic control, in \emph{Nonlinear Filtering and Stochastic Control}, Proceedings, Cortona, 1981.

\bibitem{Bismut1976} J. Bismut, Linear quadratic optimal stochastic control with random coefficients, \emph{SIAM J. Control Optim.}, {\bf 14}(3), 419-444, 1976.

\bibitem{ChenLiZhou1998}S. P. Chen, X. J. Li and X. Y. Zhou, Stochastic linear quadratic regulators with indefinite control weight costs, \emph{SIAM J. Control Optim.}, {\bf 36}(5), 1685-1702, 1998.

%\bibitem{ChenWu2022} T. Chen, Z. Wu, A general maximum principle foor partially observed mean-field stochastic system with random jumps in progressive structure, \emph{Math. Control Relat. Fields}, doi:10.3934/mcrf.2022012.

\bibitem{ChenYong2001} S. P. Chen, J. M. Yong, Stochastic linear quadratic optimal control problems, \emph{Appl. Math. Optim.}, {\bf 43}(1), 21-45, 2001.

\bibitem{ChenZhou2000} S. P. Chen, X. Y. Zhou, Stochastic linear quadratic regulators with indefinite control weight costs, \emph{II. SIAM J. Control Optim.}, {\bf 39}(4), 1065-1081, 2000.

\bibitem{ContTankov2004} R. Cont, P. Tankov, \emph{Financial Modelling with Jump Processes}, Chapmam Hall/CRC, New York, 2004.

\bibitem{Davis1977}M. H. A. Davis, \emph{Linear Estimation and Stochastic Control}, Chapman and Hall, London, 1977.

%\bibitem{Du2015}K. Du, Solvability conditions for indefinite linear quadratic optimal stochastic control problems and associated stochastic Riccati equations, \emph{SIAM J. Control Optim.}, {\bf 53}(6), 3673-3689, 2015.

%\bibitem{FOS2004} N. C. Framstad, B. \O ksendal and A. Sulem, Sufficient stochastic maximum principle for the optimal control of jump diffusions and applications to finance, \emph{J. Optim. Theory Appl.}, {\bf 121}(1), 77-98, 2004. (Errata, \emph{J. Optim. Theory Appl.}, {\bf 124}(2), 511-512, 2005.)

\bibitem{Hanson2007} F. B. Hanson, \emph{Applied Stochastic Processes and Control for Jump-Diffusions: Modeling, Analysis and Computation}, SIAM, Philadelphia, 2007.

\bibitem{HuOksendal2008} Y. Z. Hu, B. \O ksendal, Partial information linear quadratic control for jump diffusions, \emph{SIAM J. Control Optim.}, {\bf 47}(4), 1744-1761, 2008.

%\bibitem{HOP2013} S. Haadem, B. \O ksendal and F. Proske, Maximum principles for jump diffusion processes with infinite horizon, \emph{Automatica}, {\bf 49}, 2267-2275, 2013.

%\bibitem{HuangLiYong2015} J. H. Huang, X. Li and J. M. Yong, A linear-quadratic optimal control problem for mean-field stochastic differential equations in infinite horizon, \emph{Math. Control Relat. Fields}, {\bf 5}, 97-139, 2015.

\bibitem{Kou2002} S. G. Kou, A jump-diffusion model for option pricing, \emph{Manag. Sci.}, {\bf 48}(8), 1086-1101, 2002.

\bibitem{LiWuYu2018} N. Li, Z. Wu and Z. Y. Yu, Indefinite stochastic linear-quadratic optimal control problems with random jumps and related stochastic Riccati equations, \emph{Sci. China Math.}, {\bf 61}(3), 563-576, 2018.

\bibitem{LiShi2021} Z. X. Li, J. T. Shi, Linear quadratic leader-follower stochastic differential games: closed loop solvability, \emph{J. Syst. Sci. Complex.}, accepted.
	
\bibitem{LSY2016} X. Li, J. R. Sun and J. M. Yong, Mean-field stochastic linear quadratic optimal control problems: closed-loop solvability, \emph{Probab. Uncertain. Quant. Risk}, {\bf 1}(1), 24 pages, 2016.
	
\bibitem{LSY2021} X. Li, J. T. Shi and J. M. Yong, Mean-field linear-quadratic stochastic differential games in an infinite horizon, \emph{ESAIM Control Optim. Calc. Var.}, {\bf 21}, Article No. 78, 40 pages, 2021.

\bibitem{Lim2005} A. E. B. Lim, Mean-variance hedging when there are jumps, \emph{SIAM J. Control Optim.}, {\bf 44}(5), 1893-1922, 2005.

\bibitem{Lv2019} Q. Lv, Well-posedness of stochastic Riccati equations and closed-loop solvability for stochastic linear quadratic optimal control problems, \emph{J. Differential Equations}, {\bf 267}(1), 180-227, 2019.

\bibitem{Lv2020} Q. Lv, Stochastic linear quadratic optimal control problems for mean-field stochastic evolution equations, \emph{ESAIM Control Optim. Calc. Var.}, {\bf 26}, Article No. 127, 28 pages, 2020.	

%\bibitem{LvWangZhang2017}Q. Lv, T. X. Wang, and X. Zhang, Characterization of optimal feedback for stochastic linear quadratic control problems, \emph{Probab. Uncertain. Quant. Risk}, 2, 11 pages, 2017.

\bibitem{Meng2014} Q. X. Meng, General linear quadratic optimal stochastic control problem driven by a brownian motion and Poisson random martingale measure with random coefficients, \emph{Stoch. Anal. Appl.}, {\bf 32}, 88-109, 2014.

\bibitem{Merton1976} R. C. Merton, Option pricing when underlying stock returns are discontinuous, \emph{J. Finan. Econom.}, {\bf 3}(1-2), 125-144, 1976.

\bibitem{MoonChung2021} J. Moon, J. H. Chung, Indefinite linear-quadratic stochastic control problem for jump-diffusion models with random coefficients: a completion of squares approach, \emph{Mathematics}, {\bf 9}(22), Article No. 2918, 2021.

\bibitem{OS2005} B. \O ksendal, A. Sulem, \emph{Applied Stochastic Control of Jump Diffusions}, Springer, Berlin, 2005.

\bibitem{OksendalZhang2007} B. \O ksendal, T. S. Zhang, The It\^{o}-Wentzell formula and forward stochastic differential equations driven by Poisson random measures, \emph{Osaka J. Math.}, {\bf 44}, 207-230, 2007.

%\bibitem{PardouxPeng1990} E. Pardoux, S. G. Peng,  Adapted solution of a backward stochastic differential equation, \emph{Systems $\&$ Control Lett.}, {\bf 14}(1), 55-61, 1990.

%\bibitem{Peng1990} S. G. Peng, A general stochastic maximum principle for optimal control problems, \emph{SIAM J. Control Optim.}, {\bf 28}(4), 966-979, 1990.

\bibitem{Peng1992} S. G. Peng, Stochastic Hamilton-Jacobi-Bellman equations, \emph{SIAM J. Control Optim.}, {\bf 30}(2), 284-304, 1992.
		
\bibitem{Penrose1955} R. Penrose, A generalized inverse of matrices, \emph{Proc. Cambridge Philos. Soc.}, {\bf 52}, 17-19, 1955.

%\bibitem{Rishel1975} R. Rishel, A minimum principle for controlled jump processes, \emph{Lecture Notes in Economics and Mathematical Systems}, {\bf 107}, 493-508, 1975.

%\bibitem{ShenSiu2013} Y. Shen, T. K. Siu, The maximum principle for a jump-diffusion mean-field model and its application to the mean-variance problem, \emph{Nonlinear Analysis}, {\bf 86}, 58-73, 2013.

%\bibitem{ShiWu2011} J. T. Shi, Z. Wu, Relationship Between MP and DPP for the stochastic optimal control problem of jump diffusions, \emph{Appl. Math Optim.}, {\bf 63}, 151-189, 2011.

%\bibitem{Sun2017} J. R. Sun, Mean-field stochastic linear quadratic optimal control problems: open-loop solvabilities, \emph{ESAIM Control Optim. Calc. Var.}, {\bf 23}, 1099-1127, 2017.

%\bibitem{Situ1991} R. Situ, A maximum principle for optimal controls of stochastic systems with random jumps, in \emph{Proc. National Conference on Control Theory and its Applications}, Qingdao, China, 1991.
	
\bibitem{SLY2016} J. R. Sun, X. Li and J. M. Yong, Open-loop and closed-loop solvabilities for stochastic linear quadratic optimal control problems, \emph{SIAM J. Control Optim.}, {\bf 54}(5), 2274-2308, 2016.

%\bibitem{SongTangWu2020} Y. Song, S. Tang and Z. Wu, The maximum principle for progressive optimal stochastic control problems with random jumps, \emph{SIAM J. Control Optim.}, {\bf 58}£¬ 2171-2187, 2020.

\bibitem{SY2014} J. R. Sun, J. M. Yong, Linear quadratic stochastic differential games: open-loop and closed-loop saddle points, \emph{SIAM J. Control Optim.}, {\bf 52}(6), 4082-4121, 2014.
	
\bibitem{SY2018} J. R. Sun, J. M. Yong, Stochastic linear quadratic optimal control problems in infinite horizon, \emph{Appl. Math. Optim.}, {\bf 78}(1), 145-183, 2018.

\bibitem{SY2019} J. R. Sun, J. M. Yong, Linear quadratic stocahastic two-person nonzero-sum differential games: open-loop and closed-loop Nash equilibria, \emph{Stochastic Process. Appl.}, {\bf 129}(2), 381-418, 2019.
	
\bibitem{SunYong2020} J. R. Sun, J. M. Yong, \emph{Stochastic Linear-Quadratic Optimal Control Theory: Open-Loop and Closed-Loop Solutions}, Springer Briefs in Mathematics, Switzerland, 2020.

\bibitem{SunXiongYong2021} J. R. Sun, J. Xiong and J. M. Yong, Indefinite stochastic linear-quadratic optimal control problems with random coefficients: closed-loop representation of open-loop optimal controls, \emph{Anna. Appl. Probab.}, {\bf 31}(1), 460-499, 2021.

\bibitem{Tang2003} S. J. Tang, General linear quadratic optimal stochastic control problems with random coefficients: linear stochastic Hamilton systems and backward stochastic Riccati equations, \emph{SIAM J. Control Optim.}, {\bf 42}(1), 53-75, 2003.

\bibitem{Tang2015} S. J. Tang, Dynamic programming for general linear quadratic optimal stochastic control with random coefficients. \emph{SIAM J. Control Optim.}, {\bf 53}(2), 1082-1106, 2015.

\bibitem{TangHou2002} S. J. Tang, S. H. Hou, Optimal control of point processes with noisy observations: The maximum principle, \emph{Appl. Math. Optim.}, {\bf 45}, 185-212, 2002.

%\bibitem{TangLi1994} S. J. Tang, X. J. Li, Necessary conditions for optimal control of stochastic systems with random jumps, \emph{SIAM J. Control Optim.}, {\bf 32}(5), 1447¨C1475, 1994.

\bibitem{TangLiHuang2020} C. Tang, X. Q. Li and T. M. Huang, Solvability for indefinite mean-field stochastic linear quadratic optimal control with random jumps and its applications, \emph{Optim. Control Appl. Meth.}, {\bf 41}, 2320-2348, 2020.

%\bibitem{WangSunYong2019} H. X. Wang, J. R. Sun and J. M. Yong, Weak closed-loop solvability of stochastic linear-quadratic optimal control problems, \emph{Discret. Contin. Dyn. Syst. Ser. A}, {\bf 39}, 2785¨C2805, 2019.

\bibitem{Wonham1968} W. M. Wonham, On a matrix Riccati equation of stochastic control, \emph{SIAM J. Control}, {\bf 6}(1), 681-697, 1968.

\bibitem{WuWang2003} Z. Wu, X. R. Wang, FBSDE with poisson process and its application to linear quadratic stochastic optimal control problem with random jumps, \emph{Acta Automatica Sinica}, {\bf 29}(6), 821-826, 2003.

%\bibitem{Xiao2013} H. Xiao, Optimality conditions for optimal control of jump-diffusion SDEs with correlated observations noises, \emph{Mathematical Problems in Engineering}, {\bf 2013}, Article ID 613159, 7 pages, 2013.

%\bibitem{Yong2013} J. M. Yong, A linear-quadratic optimal control problem for mean-field stochastic differential equations, \emph{SIAM J. Control Optim.}, {\bf 51}(4), 2809-2838, 2013.

\bibitem{YongZhou1999} J. M. Yong, X. Y. Zhou, \emph{Stochastic Controls: Hamiltonian Systems and HJB Equations}, Springer-Verlag, New York, 1999.

\bibitem{Yu2013}Z. Y. Yu, Equivalent cost functionals and stochastic linear quadratic optimal control problems, \emph{ESAIM Control Optim. Calc. Var.}, {\bf 19}, 78-90, 2013.

\bibitem{ZhangDongMeng2020} F. Zhang, Y. C. Dong and Q. X. Meng, Backward stochatic Riccati equation with jumps associated with stochastic linear quadratic optimal control with jumps and random coefficients, \emph{SIAM J. Control Optim.}, {\bf 58}(1), 393-424, 2020.

\end{thebibliography}
\end{document}